\documentclass[12pt]{amsart}      
\usepackage{amsmath, amstext, amsbsy, amssymb}

\newtheorem{theorem}{Theorem}[section]
\newtheorem{lemma}[theorem]{Lemma}

\newtheorem{proposition}[theorem]{Proposition} 
\newtheorem{corollary}[theorem]{Corollary}  

\theoremstyle{definition}

\theoremstyle{remark}
\newtheorem{remark}[theorem]{Remark}

\numberwithin{equation}{section}

\newcommand{\ch}{\mbox{ch}}

\newcommand{\tFG}{{\mathcal F}_{\G}^-}
\newcommand{\FGG}{\overline{\mathcal F}_{\G} }
\newcommand{\tFGG}{\overline{\mathcal F}_{\G}}

\newcommand{\G}{\Gamma}
\newcommand{\Gbar}{\overline{\Gamma}^*}
\newcommand{\Gm}{{\Gamma}_m}
\newcommand{\Gn}{{\Gamma}_n}

\newcommand{\gs}{\mathfrak g}
\newcommand{\g}{\gamma}
\newcommand{\ga}{\gamma}
\newcommand{\thg}{\widehat{\mathfrak h}_{\G, \wt}[-1]}
\newcommand{\hg}{\widehat{\mathfrak h}_{\G, \wt}}

\newcommand{\thhg}{\widehat{\widehat{\mathfrak g}}[-1]}
\newcommand{\la}{\lambda}
\newcommand{\loopg}{\widehat{\mathfrak g}}
\newcommand{\tloopg}{\widehat{\mathfrak g}[-1]}

\newcommand{\be}{\beta}

\newcommand{\ep}{\epsilon}
\newcommand{\RG}{R_{\G}}
\newcommand{\tRG}{R^-_{\G}}
\newcommand{\RGG}{\overline{R}_{ \G}}
\newcommand{\tRGG}{\overline{R}_{ \G}}
\newcommand{\RGO}{R^0_{ \G}}
\newcommand{\Rz}{R_{\mathbb Z} (\G)}
\newcommand{\Rzz}{\overline{R}_{\mathbb Z} (\G)}
\newcommand{\SG}{ S_\G }
\newcommand{\tSG}{ S^-_\G }

\newcommand{\tSGG}{\overline{S}_\G }
\newcommand{\SGO}{ S^0_\G }

\newcommand{\tVG}{V_{\G}^-}

\newcommand{\tVGG}{\overline{V}_{ \G}}

\newcommand{\wt}{\xi}
\newcommand{\Z}{\mathbb Z}

\def\tS#1{\widetilde{S}_{#1}}
\def\tG#1{\widetilde{\Gamma}_{#1}}
\newcommand{\htimes}{\hat{\times}}
\newcommand{\Rtz}{R_{\mathbb F_2} (\G)}
\newcommand{\Rtzbar}{R_{\mathbb F_2} (\Gbar)}
\newcommand{\Rtzh}{\hat{R}_{\mathbb F_2}^-(\G)}
\newcommand{\Rtzhbar}{\hat{R}_{\mathbb F_2}^-(\Gbar)}
\newcommand{\oa}{\overline{\alpha}}
\newcommand{\ob}{\overline{\beta}}

\begin{document}

\title[Twisted vertex operators and McKay correspondence]
    {Twisted vertex representations via spin groups and 
the McKay correspondence}
\author{Igor B. Frenkel}
\address{Department of Mathematics,
   Yale University, 
   New Haven, CT 06520}
\author{Naihuan Jing}
\address{Department of Mathematics,
   North Carolina State Univer\-sity,
   Ra\-leigh, NC 27695-8205}
\email{jing@math.ncsu.edu}
\thanks{Research of Frenkel is supported by
NSF grant DMS-9700765;
research of Jing is supported by NSF grant 
DMS-9970493.}
\author{Weiqiang Wang}
\address{Department of Mathematics,
   North Carolina State Univer\-sity,
   Ra\-leigh, NC 27695-8205}
\email{wqwang@math.ncsu.edu} 
\begin{abstract}
We establish a twisted analog of our recent work on 
vertex representations and the McKay correspondence. For
each finite group $\Gamma$ and a virtual character of
$\Gamma$ we construct twisted vertex operators on
the Fock
space spanned by the super spin characters of the spin wreath products
$\Gamma\wr\widetilde{S}_n$ of $\Gamma$ and a double cover
of the symmetric group $S_n$ for all $n$. 
When $\Gamma$ is a subgroup of $SL_2(\mathbb C)$ with the McKay
virtual character,
our construction gives a group theoretic
realization of the basic representations of the twisted
affine and twisted toroidal algebras.
When $\Gamma$ is an arbitrary finite group and the virtual character
is trivial, our vertex operator
construction yields the spin character tables for
$\Gamma\wr\widetilde{S}_n$.
\end{abstract}

\maketitle

\section{Introduction} \label{S:intro}

The connection among the direct sum of Grothendieck groups of the
symmetric groups $S_n$ for all $n$ and the theory of
symmetric functions \cite{M, Z} has
a simple interpretation in terms of
a Heisenberg algebra and
vertex operators (\cite{F1}, see part one of \cite{J}). 
In the recent works \cite{W, FJW1} we have
realized a generalization of such a connection by substituting the symmetric group $S_n$ with the wreath product $\Gn = \G \wr S_n$
associated to an arbitrary finite group $\G$. 
Moreover, we introduced a crucial modification of this connection that, in the case when $\G$ is a finite
subgroup of $SL_2(\mathbb C)$, yields a 
group theoretic realization of the affine Lie algebra
$\loopg$ \cite{FK, Se} and of the toroidal Lie algebra
$\widehat{\loopg}$ \cite{F2, MRY}, where $\mathfrak g$ is a complex
simple Lie algebra of ADE type whose Dynkin diagram is related to
$\G$ via the McKay correspondence \cite{Mc}.

The main goal of the present work is to extend the above results 
to realize the twisted basic representation of an affine Lie
algebra $\tloopg$ and its toroidal counterpart by means of a spin cover
$\tG n$ of the wreath product $\Gn$ associated to a subgroup $\G$ of $SL_2(\mathbb C)$.

The twisting of the basic representation of the affine Lie algebra
under consideration is determined by the multiplication by $-1$
on the Cartan subalgebra $\mathfrak h$ of $\mathfrak g$, and  can
be viewed as an odd counterpart of the even (untwisted) case. This
twisting was originally introduced as the first step towards the
construction of the Moonshine module for the Monster group in
\cite{FLM1, FLM2}. As in the homogeneous case one starts with a
representation of the Heisenberg subalgebra $\widehat{\mathfrak
h}[-1]$ and reconstructs the rest of the twisted affine Lie
algebra $\tloopg$ using the twisted vertex operators.

The representation theory of the spin group $\tS n$ which is a
double cover of the symmetric group $S_n$ was initiated by 
I. Schur \cite{Sc} (also see \cite{Jo} for an
exposition). Its connection with vertex operators was further
studied in \cite{J}. These
results will play an important role in our present work. 
The representation theory of $\tG n$ was
also studied in \cite{HH} from a Hopf algebra viewpoint.

In order to work effectively only with the spin
representations of $\tG n$, i.e, those which do not factor through
$\Gn $, we adopt the approach of \cite{Jo} by introducing a
superalgebra structure on the group algebra of $\tG n$ and
consider its supermodules. It turns out that the superstructure is
preserved under the main operations such as induction and
restriction. The direct sum of the Grothendieck groups
of spin supermodules of $\tG n$ carries a natural Hopf algebra
and we remark that a Hopf algebra was constructed in \cite{HH} on a different space.
 This allows us to realize the vertex operators acting
in the twisted vertex representations  constructed from the
sum of the Grothendieck rings. Our group theoretic 
method naturally recovers the basic representations of twisted
affine Lie algebras $\tloopg$ \cite{LW, FLM1, FLM2}. As in \cite{FJW1}
we realize this by introducing a modified bilinear form 
associated to the McKay virtual character $\xi$ which is twice the
trivial character minus the character of the two-dimensional
natural representation of $\G$ in $SL_2 (\mathbb C)$. 

Much of our construction is valid for an arbitrary finite group
$\G$ and we have introduced the modified bilinear form associated
to an arbitrary virtual character $\xi$ of $\G$ as well. In the
special case when $\xi$ is the trivial character the twisted
vertex operators generate an infinite dimensional generalized Clifford
algebra, which recovers the twisted boson-fermion correspondence. We further
obtain the super character tables of the spin group $\tG n$ for all
$n$, generalizing the results of \cite{J}.

One may generalize the results of this paper to the quantum case
as it was done in \cite{FJW2} for the homogeneous picture of quantum affine
algebras
\cite{FJ}. Our results 
also suggest that various previous constructions
associated to (quantum) vertex representations admit remarkable
interpretation via Grothendieck rings of certain finite groups
which are variations of wreath products, though every new step in
this direction is unpredictable and brings new surprises. 
It is a
very interesting and challenging problem to find such a group
theoretic realization.

The organization of the paper is as follows. In Sect. 2 we present the representation theory and structures of the
spin group $\tG n$.
In Sect. 3 we review superalgebras and supermodules
and define the Hopf algebra of the super
spin characters of $\tG n$. In Sect. 4 we introduce the
weighted bilinear forms in the Grothendieck rings of supermodules and construct basic spin supermodules.
In Sect. 5 we define the twisted Heisenberg algebras and their
Fock spaces. In Sect. 6 we establish the isometry between the sum
of Grothendieck rings of supermodules of $\tG n$
and the Fock space of a twisted Heisenberg algebra. In
Sect. 7 we construct twisted vertex operators via the induction
and restriction functors on the Grothendieck rings. In Sect. 8 we
obtain the twisted basic representation of the affine Lie algebras
$\tloopg$ and the corresponding toroidal algebras. In Sect. 9 we
derive the super spin character tables of $\tG n$ for all $n$
from the twisted boson-fermion correspondence.

\section{A double cover of the wreath product} \label{sect_wreath}
\subsection{The spin group $\tS n$}
In this subsection we discuss some of the basic properties of the 
double covers of the symmetric group, which were
introduced by Schur in his seminal paper \cite{Sc}. We
will adopt the modern account \cite{Jo} of Schur's theory.

Let $S_n$ be the symmetric group of $n$ letters, and we use
the convention of multiplying permutations from right to left
(different from \cite{Sc, Jo}). The spin group $\tS n$ 
is the finite group generated by $z$ and 
$t_i, i=1, \cdots, n-1$ subject to the relations:
\begin{align}\label{E:defrel1}
&z^2=1,  \quad t_i^2=(t_it_{i+1})^3=z, \\ \label{E:defrel2}
&t_it_j=zt_jt_i, \qquad i>j+1, \\ \label{E:defrel3}
&zt_i=t_iz.
\end{align}

Let $\theta_n$ be the homomorphism from $\tS n$ to $S_n$ sending $t_i$ to 
the transposition $(i, i+1)$ and $z$ to $1$. 
We see that $\tS n$ is a central extension of $S_n$ by
the cyclic group $\mathbb Z_2$:
\begin{equation*}
1 \longrightarrow\mathbb Z_2 \stackrel{\iota}{\longrightarrow}
\tS n \stackrel{\theta_n}{\longrightarrow} S_n \longrightarrow 1,
\end{equation*}
where the embedding $\iota$ sends the order $2$ element in
$\mathbb Z_2$ to $z$.
Schur \cite{Sc} determined that
$H^2(S_n, \mathbb C^*)\simeq \mathbb Z_2$ for $n>3$. The group $\tS n$ is one of the two double covers of the symmetric group $S_n$ ($n>3$). Our results in this paper can be easily translated to the other double cover
(cf. \cite{Sc, J}).

The group $\tS n$ has a parity given as follows.
Let $d$ be the homomorphism from the free group generated by
$\{t_i, z\} (i=1, \cdots, n-1) $ to $\mathbb Z_2$ by
$d(t_i)=1$, $i=1, \cdots, n-1$ and $d(z)=0$.
It is easily seen that $d$ preserves the relations 
(\ref{E:defrel1}-\ref{E:defrel3}). Thus it defines a homomorphism
from $\tS n$ to $\mathbb Z_2$,
which we still denote by $d$. 
An element $x\in\tS n$ is called {\it even} (resp. {\it odd})
if $d(x)=0$ (resp. $d(x)=1$). The parity in $\tS n$ given by $d$ lifts the 
standard notion of even and odd permutations in the symmetric 
group $S_n$.

The spin group $\tS n$ has a cycle presentation due to J.~H.~Conway
and others (see \cite{Ws}). Embed $\tS n$ into
$\tS{n+1}$ by identifying their first $n-1$ 
generators $t_i, i=1, \cdots, n-1$. 
For $i=1, \cdots, n$ we define $x_i=t_it_{i+1}
\cdots t_n\cdots t_{i+1}t_i \in \tS{n+1}$. For a sequence
$i_1, \cdots, i_m$ of distinct integers
from $\{1, 2, \cdots, n\}$ we can define 
cycles in $\tS{n}$ as follows.
\begin{equation}
[i_1i_2\cdots i_m]=\begin{cases} z, & m=1,\\
   x_{i_1}x_{i_m}x_{i_{m-1}}\cdots x_{i_1}, & 1<m\le n.
\end{cases}
\end{equation}
It is known that $\theta_n([i_1i_2\cdots i_m])=
(i_1i_2\cdots i_m)$ and $\theta_{n+1}(x_i)=(i, n+1)$. 
We list some useful identities for the cycles.
\begin{align} \label{E:commut1}
&x_j[i_1i_2\cdots i_m]=z^{m-1}[i_1i_2\cdots i_m]x_j,\quad (j\neq i_s),
\quad x_j^2=z,\\ 
&[i_1i_2\cdots i_m]^{-1}=[i_m\cdots i_2i_1],\\
&[i_1i_2\cdots i_m]=z^{m-1}[i_2i_3\cdots i_mi_1], \label{E:commut2}\\
&[i_1i_2\cdots i_m][j_1j_2\cdots j_k]=z^{(m-1)(k-1)}
[j_1j_2\cdots j_k][i_1i_2\cdots i_m], \label{E:commut}\\
&[i,i+1,\cdots, i+j-1]=z^{j-1}t_it_{i+1}\cdots t_{i+j-2},
\end{align}
where the cycles $[i_1i_2\cdots i_m]$ and $[j_1j_2\cdots j_k]$ are disjoint.

\begin{proposition} \cite{Jo}
\label{P:present}
Each element of
$\tS n$ can be presented as
$$ z^p[i_1i_2\cdots i_m][j_1j_2\cdots j_k]\cdots ,$$
where $\{i_1\cdots i_m\}, \{j_1\cdots j_k\}, \cdots $
is a partition of the set $\{1, 2, \cdots, n\}$ and $p=0, 1$. If $z^pc_1c_2\cdots c_l=z^{p'}c'_1 c'_2\cdots c'_{l'}$
are two expressions of the same element in terms of
cycles $c_i$ and $c_i'$, then
$l=l'$ and there is a permutation $\sigma\in S_l$ such that
$$ c_i'=c_{\sigma(i)} z^{m_i}, \qquad m_i\equiv |c_i|-1 (mod\, 2),
$$
where $|c_i|$ denotes the length of the cycle $c_i$. 
Moreover if $[i_1\cdots i_k]=z^m[j_1\cdots j_k]$, then
$j_s=\sigma(i_s)$ for a cyclic permutation $\sigma$
of $\{i_1, \cdots, i_k\}$.
\end{proposition}

Let $\lambda$ be a partition and we identify $\la$ with its
Young diagram consisting of $l$ rows of 
$\la_1, \cdots, \la_l$ squares respectively aligned to the left. 
A {\it tableau} $T_{\la}$ of shape 
$\lambda$ is a numbering of the squares with integers
$1, 2, \cdots, |\la|$, each appearing exactly once. 
For each tableau $T_{\la}$ of shape $\la$ with a numbering
$a_{11}, \cdots a_{1\la_1}$, $a_{21}, \cdots, a_{2\la_2}$,
$\cdots, a_{l1}, \cdots, a_{l\la_l}$ we define 
the element $t_{\la}$ of $\tS n$ to be
\begin{equation}\label{E:tableau1}
t_{\la}=[a_{11}\cdots a_{1\la_1}][a_{21}\cdots a_{2\la_2}]\cdots
[a_{l1}\cdots a_{l\la_l}].
\end{equation}
The permutation $\prod_{i=1}^{l}(a_{i1}\cdots a_{i\la_i})$
associated with $t_{\la}$
will be denoted by $s(\la)$.
It follows from Proposition \ref{P:present} that the general
element in $\tS n$ is of the form $z^pt_{\la}$. For a permutation
$s\in S_n$ we also define
$t_{\la}^{s}=\prod_{i=1}^l[s(a_{i1})\cdots s(a_{i\la_i})]$.  

The following can be checked by
induction using (\ref{E:commut1}) and (\ref{E:commut}).

\begin{lemma} \label{L:conj}
 For any two elements  $t_{\la}, t_{\mu}$ in
$\tS n$ associated to tableaux $T_{\la}$ and $T_{\mu}$
we have that
\[ t_{\mu}t_{\la}t_{\mu}^{-1}=z^{d(\la)d(\mu)}
t_{\la}^{s(\mu)}. 
\]
\end{lemma}

\subsection{The spin wreath product $\tG n$} 
In this subsection we introduce the main finite group
$\tG n$ in this work, and extend our discussion from 
$\tS n$ to $\tG n$.

  Let $\Gamma$ be a finite group with $r +1$ conjugacy classes.
We denote by $\Gamma^*=\{\g_i\}_{i=0}^{r}$ the set of complex 
irreducible characters, where
$\g_0$ stands for the trivial character, and by 
$\Gamma_*$ the set of 
conjugacy classes. The character value $ \gamma(c) $
of $\g \in \G^*$ at a conjugacy class $c\in \Gamma_*$
yields the character table $ \{ \gamma(c) \} $ of $\G$.

Let $R(\G)=\oplus_{i=0}^{r} \mathbb C\g_i$ be the 
space of complex valued class functions on $\G$. 
For $c \in \G_*$
let $\zeta_c$ be the order of the centralizer of an element
in the class $c$, so the
order of the class is then $|\G |/\zeta_c$. 
The usual bilinear form on $R(\G )$ is defined
as follows:

\begin{eqnarray*}
\langle f, g \rangle_{\G} = \frac1{ | \G |}\sum_{x \in \Gamma}
          f(x) g(x^{ -1})
 = \sum_{c \in \Gamma_*} \zeta_c^{ -1} f(c) g(c^{ -1}),
\end{eqnarray*}
where $c^{ -1}$ denotes the conjugacy class 
$\{ x^{ -1}, x \in c \}$. Clearly $\zeta_c = \zeta_{c^{-1}}$.
We will often write $\langle \ , \ \rangle$
for $\langle \ , \ \rangle_{\G }$ when no ambiguity may arise. 
It is well known that 
\begin{eqnarray}
  \langle \g_i, \g_j \rangle &= & \delta_{ij} , \nonumber \\
  \sum_{ \g \in \G^*} \g (c ')  \g ( c^{ -1})
    &= & \delta_{c, c '} \zeta_c, \quad c, c ' \in \G_*.  \label{eq_orth}
\end{eqnarray}
Thus $\Rz=\oplus_{i=0}^{r} \mathbb Z\g_i$ 
endowed with this bilinear form
becomes an integral lattice in $R( \G )$.

Given a positive integer $n$, let $\Gamma^n = \Gamma \times \cdots
\times \Gamma$ be the $n$-th direct product of $\Gamma$, and let
$\G^0 $ be the trivial group. 
The spin group $\tS n$ acts on 
$\Gamma^n$ through the action of the group
$S_n$ by permuting the indices: 
$t_{\la} (g_1, \cdots, g_n)
  = (g_{s(\la)^{-1} (1)}, \cdots, g_{s(\la)^{ -1} (n)}),
$ and $z(g_1, \cdots, g_n)=(g_1, \cdots, g_n)$.
The wreath product $\tG n=\G\wr \tS n$ of $\Gamma$ with $\tS n$ is defined to be
the semi-direct product
$$
 \tG n =\G^n\rtimes \tS n= \{(g, t) | g=(g_1, \cdots, g_n)\in {\Gamma}^n,
t\in \tS n \}
$$
 with the multiplication
$$
(g, t)\cdot (h, s)=(g \, {t} (h), ts) .
$$
Note that $\tG n$ reduces to $\tS n$ when $\G$
is trivial. Clearly $\tG n$ is a central
extension of $\Gn$ by $\mathbb Z_2$ and $|\tG n|=2n!|\G|^n$.

We define a a parity for $\tG n$ 
by extending the parity of $\tS n$. 
Let $d: \tG n\longrightarrow \mathbb Z_2=\{0, 1\}$ be the homomorphism from $\tG n$ to $\mathbb Z_2$ given by
\begin{equation}\label{E:parity}
d(g, t_i)=1, \qquad d(g, z)=0.
\end{equation}
Clearly the degree $0$
subset $\tG n^0$ is the wreath product $\G\wr \widetilde{A}_n$, where
$A_n$ is the alternating group, and the 
degree one part $\tG n^1$ is the complementary subset.

Let $\tau$ be a section from $\Gn$ to $\tG n$ such
that $\theta\tau=1$. 
An element $x\in \Gn$ is called {\it split} 
if $\tau(x)$
is not conjugate to $z\tau(x)$. Otherwise $x$ is said to be
{\it non-split}. Clearly this definition does not depend
on the choice of the section $\tau$ and
two conjugate elements 
are simultaneously split or non-split. 
A conjugacy class of $\Gn$ is called
{\it split} if its
elements are split.  We will also
say that an element $x\in\tG n$ is split (resp. non-split) 
if $\theta(x)$ is split (resp. non-split). 
Clearly the  class $C_{\rho}$ splits if and only if the preimage $\theta_n^{-1}(C_{\rho})$ splits into two conjugacy 
classes in $\tG n$.

Any representation $\pi$ of $\Gn$ can be viewed as
a representation of $\tG n$. Such a representation
$\pi$ of $\tG n$ satisfies the property 
$\pi(z)=Id$.  A  representation
$\pi$ of $\tG n$ is called {\it spin} if 
$\pi(z)=-Id.$
It follows that the characters of spin representations vanish on non-split classes. In this paper
we only consider spin representations.

We remark that spin representations are sometimes referred as
negative or projective
representations in the literature.

\subsection{Conjugacy classes of $\tG n$} \label{S:conjclass} 
Let $\la=(\la_1, \la_2, \cdots, \la_l)$ be a partition
of the integer $|\la|=\la_1+\cdots+\la_l$, where 
$\la_1\geq \dots \geq \la_l \geq 1$.
The integer $l$ is called the {\em length} of the partition
$\la $ and is denoted by $l (\la )$. 
We will identify the partition $(\la_1, \la_2, \cdots, \la_l)$ with
$(\la_1, \la_2, \cdots, \la_l, 0, \cdots, 0)$.
We will also make use of another notation for partitions:
$$
\la=(1^{m_1}2^{m_2}\cdots) ,
$$
where $m_i$ is the number of parts in $\la$ equal to $i$.
The number $n(\la')$ is defined to be $\sum_i\binom{\la_i}2$, where $\la'$ is the dual partition
associated to $\la$. 
We will use the dominance order on partitions. For two partitions $\la$ and $\mu$ we write $\la\geqslant\mu$  
if $\la_1\geq \mu_1$, $\la_1+\la_2\geq \mu_1+\mu_2$, etc.

A partition $\la$ is called {\it strict} if its parts are distinct
integers (excluding the trivial parts of zero), in which case
all the multiplicities $m_i$ are $1$.

We will use partitions indexed by $\G_*$ and $\G^*$. For
a finite set $X$ and $\rho=(\rho(x))_{x\in X}$ a family
of partitions indexed by $X$, we write 
$$\|\rho\|=\sum_{x\in X}|\rho(x)|.$$
It is convenient to regard $\rho=(\rho(x))_{x\in X}$
as a partition-valued function on $X$. 
We denote by $\mathcal{P}(X)$ the set of all partitions indexed by $X$ 
and by $\mathcal{P}_n(X)$ the set of all partitions
in $\mathcal{P}(X)$ such that $\|\rho\|=n$. The total
number of parts, denoted by $l(\rho)=\sum_xl(\rho(x))$, in the partition-valued function $\rho=(\rho(x))_{x\in X}$
is 
called the length of $\rho$.
The {\it dominance order} is extended to partition-valued functions as follows. We define
$\rho\geqslant \pi$ if $\rho(x)\geqslant \pi(x)$
for each $x$. We say that
$\rho\gg\pi$ if $\rho(x)\geqslant \pi(x)$ and $\rho(x)\neq\pi(x)$ for {\it each} $x\in X$.
For a partition-valued function $\rho$ we define
\begin{equation}\label{E:n}
n(\rho')=\sum_cn(\rho(c)')=\sum_{c, i}\binom{\rho_i(c)}2.
\end{equation}

Let $\mathcal{OP}(X)$ be the set of partition-valued functions $(\rho(x))_{x\in X}$ in
$\mathcal{P}(X)$ such that all parts of the partitions
$\rho(x)$
are odd integers,
and let $\mathcal{SP}(X)$ be the set of partition-valued
functions $\rho: X\longrightarrow \mathcal P$ such that
each partition $\rho(x)$ is strict. When $X$ consists
of a single element, we will omit $X$ and simply write
$\mathcal P$ for $\mathcal P(X)$, thus the notation
$\mathcal{OP}$ or $\mathcal{SP}$ will be used similarly.   

\begin{lemma} \label{L:euler}
$|\mathcal{OP}_n(X)|=|\mathcal{SP}_n(X)|$.
\end{lemma}
\begin{proof} The generating function of the cardinalities of strict partition-valued functions is
\begin{equation*}
\prod_{x\in X}\prod_{n=1}^{\infty}\frac 1{1-q^{2n-1}_x}
=\prod_{x\in X}\prod_{n=1}^{\infty}
\frac{1-q^{2n}_x}{(1-q^{2n-1}_x)(1-q^{2n}_x)}
=\prod_{x\in X}\prod_{n=1}^{\infty}(1+q^n_x),
\end{equation*}
which is the generating function of $\mathcal{SP}_n(X)$.
\end{proof}

We also define a parity on partitions. 
For each partition $\la$ we define
$d(\la)=|\la|-l(\la)$. For a partition-valued function
$\rho=(\rho(x))_{x\in X}$ we define 
$d(\rho)=\sum_x |\rho(x)|=\|\rho\|-l(\rho)$.
It is clear that the conjugacy class
of type $\la$ in $S_n$ is even if and only if $d(\la)$
is even. 
We define the parity of the partition-valued
function $\rho$ to be the parity of $d(\rho)$.
We define
\begin{align}\label{E:parity1}
\mathcal P_n^0(X)&=\{\la\in \mathcal P_n(X)|\quad d(\rho)\equiv 0 (mod\, 2)\},\\ \label{E:parity11}
\mathcal P_n^1(X) &=\{\la\in \mathcal P_n(X)|\quad d(\rho)\equiv 1 (mod\, 2)\},
\end{align}
and define $\mathcal {SP}_n^i(X)=\mathcal {P}_n^i(X)\cap
\mathcal{SP}_n(X)$ for $i=0, 1$.

We now recall the description of conjugacy classes of ${\Gamma}_n$ \cite{M}. Let $x=(g, \sigma)$ be an
element in a conjugacy class of $\Gn$, where
$g=(g_1, \cdots, g_n)$. For each 
cycle $y=(i_1 i_2 \cdots i_k)$ in the permutation
$\sigma$
the element $g_y=g_{i_k} g_{i_{k -1}}
\cdots g_{i_1} \in \Gamma$ is called the {\em cycle-product}
of $x$ corresponding to the cycle $y$. For each $c\in\G_*$
and $i\geq 0$ let
$m_i(c)$ be the number of $i$-cycles in the permutation
$\sigma$ such that
the cycle products $g_y$ lie in the conjugacy class
$c$. Then $c\rightarrow
\rho(c)=(1^{m_1(c)}2^{m_2(c)}\cdots)$ defines 
a partition-valued function on $\G_*$. It is known that 
the partition-valued function $(\rho(c))_{c\in\G_*}$
is in one-to-one correspondence to the conjugacy class of $x=(g, \sigma)$ in $\Gn$
and is called the {\it type} of the conjugacy class. 
We will also say that an element has conjugacy type $\rho$ if
this element is contained in the conjugacy class.
 
Let $(-1)^d$ be the representation of
$\tG n$ given by $x\longmapsto (-1)^{d(x)}$.
A representation $\pi$ of $\tG n$ is called a {\it double spin} representation
if 
\[ (-1)^{d}\pi \simeq\pi.
\]
If $\pi'=(-1)^{d}\pi \ne\pi$, then $\pi'$ and $\pi$ are called
{\it associate spin} representations of $\tG n$.

The following result was proved in \cite{J} for a  
double cover of any finite group.

\begin{proposition} The number of split conjugacy classes of
$\Gn$
is equal to the number of irreducible spin representations of $\tG n$.
\end{proposition}

\subsection{Split conjugacy classes of $\tG n$}\label{S:splitclass}
We fix an order of conjugacy classes of $\G$:
$c^0=\{1\}, c^1, \cdots, c^r$.
For each partition-valued function 
$\rho=(\rho(c))\in\mathcal{P}_n(\G_*)$, we 
let $t_{\rho(c^i)}$ be the element
of $\tS n$ associated to a tableau $T_{\rho(c^i)}$
of shape $\rho(c^i)$ using the numbers
$\sum_{j\leq i-1}|\rho(c^j)|+1, \cdots, 
\sum_{j\leq i}|\rho(c^j)|$ and we define
the element $t_{\rho}$ to be
\begin{equation}\label{E:specelt3}
t_{{\rho}}=t_{\rho(c^0)}t_{\rho(c^1)}\cdots t_{\rho(c^r)},
\end{equation}
which depends on the sequence $T_{\rho}$
of the tableaux
$T_{\rho(c^0)}, \cdots, T_{\rho(c^r)}$. We remark that
the general element of $\tG n$ is of the form
$(g, z^{p}t_{\rho})$, where $\rho$ is the type of the
conjugacy class of $(g, z^{p}t_{\rho})$. 

The following theorem is well-known in the case
of $\G=\{1\}$ (cf. \cite{Jo} and \cite{St}).

\begin{theorem} \label{T:class}
Let $\rho=(\rho(c))_{c\in\G_*}$ be the type of a conjugacy class 
$C_{\rho}$ in $\Gn$. Then the preimage
$\theta_n^{-1}(C_{\rho})$ splits into two conjugacy classes
in $\tG n$ if and only if 

(1) when the class $C_{\rho}$ is even and all the $\rho(c)$ $(c\in\G_*)$ are 
partitions with odd integer parts, i.e., 
$\rho\in \mathcal{OP}_n(\G_*)$;

(2) when the class $C_{\rho}$ is odd and all the $\rho(c)$
$(c\in\G_*)$ are strict
partitions, i.e., $\rho\in \mathcal{SP}_n^1(\G_*)$.
\end{theorem}
\begin{proof}
(1) Let $d(\rho)$ be even and let each
partition $\rho(c)$ have odd integer parts. 
Assume on the contrary that $(g, t_{\rho})$ and $z(g, t_{\rho})$ are conjugate in
$\theta^{-1}_n(C_{\rho})$, where $t_{\rho}$
is associated to a sequence of
tableaux (see
Eqn. (\ref{E:specelt3})). 
Then for some $(h, t_{\mu})\in\tG n$
\begin{align}\nonumber
&(h, t_{\mu})(g, t_{\rho})(h, t_{\mu})^{-1}\\ 
=&(h\cdot s(\mu)(g)\cdot s(\mu)s(\rho)s(\mu)^{-1}(h^{-1}), t_{\mu}t_{\rho}
t_{\mu}^{-1}) \nonumber\\ \label{E:split1}
=&(h\cdot s(\mu)(g)\cdot s(\rho)(h^{-1}), t_{\mu}t_{\rho}
t_{\mu}^{-1})=(g, zt_{\rho}),
\end{align}
where we have used the fact that $s(\rho)s(\mu)=s(\mu)s(\rho)$.
It follows from Lemma \ref{L:conj} that
$zt_{\rho}=z^{d(\rho)d(\mu)}t_{\rho}^{s(\mu)}=
t_{\rho}^{s(\mu)}$,
since $d(\rho)=0 \, (mod\, 2)$. 
Let $t_{\rho}=c_1c_2\cdots c_l$ and $t_{\rho}^{s(\mu)}=
c_1'c_2'\cdots c_l'$ be their cycle
representations. Then $c_i=z^{m_i}c_{\nu(i)}'$ and $m_i=|c_i|-1
 (mod \, 2)$ for some $\nu\in S_l$
by Proposition \ref{P:present}. Since each cycle length
$|c_i|$ is odd, all the cycles mutually commute with each other.
Substituting $c_i=z^{m_i}c_{\nu(i)}'$ back and  
rearranging the cycles, we have 
\begin{equation*}
1=z^{1+\sum_i{(|c_i|-1})}=z^{1+d(\rho)}=z,
\end{equation*}
which is a contradiction.

Now suppose that for some $c\in \G_*$
there is an even cycle in $\rho(c)$ of the class
$C_{\rho}$ of type $\rho$. That is, there is an element $(g, t_{\rho})
\in \theta^{-1}(C_{\rho})$
such that $t_{\rho}=\cdots [i_1i_2\cdots i_{2k}]\cdots$.
Consider the element $(h, t_{\mu})\in \tG n$, 
where $t_{\mu}=[i_1i_2\cdots i_{2k}]$
and $h=(h_1, \cdots, h_n)$ with $h_j=1$ for $j\neq i_s$
and $h_{i_s}=g_{i_s}, s=1, \cdots, 2k$. We claim that
\begin{equation*}
(h, t_{\mu})(g, t_{\rho})(h, t_{\mu})^{-1}
=(hs(\mu)(g)s(\rho)(h^{-1}), t_{\mu}t_{\rho}
t_{\mu}^{-1})=(g, zt_{\rho}),
\end{equation*}
which is shown by two steps. First we consider the $j$th component
of $hs(\mu)(g)s(\rho)(h^{-1})$ in $\G^n$. It equals 
$1\cdot g_j\cdot 1=g_j$ when $j\neq i_s$, and it 
equals 
$g_{i_s}g_{i_{s-1}}g_{i_{s-1}}^{-1}=g_{i_s}$ for $j=i_s$. 
Secondly we have
\begin{equation*}
t_{\mu}t_{\rho}
t_{\mu}^{-1}=\cdots [i_2i_3\cdots i_{2k}i_1]\cdots
=zt_{\rho}
\end{equation*}
by using Eqn. (\ref{E:commut2}) and $d(\rho)\equiv 0\, (mod\, 2)$ again. Thus
$(g, t_{\rho})$ is conjugate to $z(g, t_{\rho})$. Therefore all
partitions $\rho(c)$ must be from $\mathcal{OP}(\G_*)$
if $\theta^{-1}(C_{\rho})$ splits.

(2) Let $d(\rho)$ be odd. Assume all partitions $\rho(c)$
are strict partitions. If on the contrary 
$(g, t_{\rho})$ is conjugate to $z(g, t_{\rho})$, then using
$d(\rho)=1$ we have as in
(\ref{E:split1}) that $zt_{\rho}=z^{d(\mu)}t_{\rho}^{s(\mu)}$
for the permutation $s(\mu)\in S_n$
associated to some $\mu\in\mathcal{P}_n(\G_*)$. Let
$t_{\rho}=c_1c_2\cdots c_l$ and $t_{\rho}^{s(\mu)}=
c_1'c_2'\cdots c_l'$ be their cycle
representations. Then $c_i=z^{|c_i|-1}c_{\nu(i)}'$ for
some $\nu\in S_l$. Since $s(\mu)$ cyclically permutes
the indices in each cycles of $s(\rho)$ we have $d(\mu)=d(\rho)$. 
On the other hand, note that each cycle $c_i$ corresponds
to one part in $\rho(c)$ for some $c\in\G_*$ and
any conjugation of $c_i$ still corresponds to a part in
the same $\rho(c)$. When we plug the equations
$c_i=z^{|c_i|-1}c_{\nu(i)}'$ back to $zt_{\rho}=z^{d(\mu)}t_{\rho}^{s(\mu)}$
we see that $\nu$ is actually the identity since
$\rho(c)$ is strict. Therefore 
$z^{1+\sum_i(|c_i|-1)}=z^{d(\mu)}$. Then $d(\rho)=\sum_i(|c_i|-1)\equiv
d(\mu)+1 \, (mod\, 2)$, which is a contradiction.
Hence 
$\theta^{-1}_n(C_{\rho})$ splits.

Now suppose $\theta^{-1}_n(C_{\rho})$ splits.
If there are two identical parts in $\rho(c)$ for some conjugacy class
$c\in \G_*$, say $t_{\rho}=
\cdots [i_1\cdots i_k][j_1\cdots j_k]\cdots$
for $(g, t_{\rho})\in \tG n$. Then the cycle-products of these
two identical parts are conjugate, i.e., 
there exists an element $x\in\G$ such that
\begin{equation}\label{E:conj-pf1}
xg_{j_k}g_{j_{k-1}}\cdots g_{j_1}x^{-1}
=g_{i_k}g_{i_{k-1}}\cdots g_{i_1}.
\end{equation}
Consider the element
$(h, t_{\mu})$ such that $t_{\mu}=[i_1j_1]\cdots [i_kj_k]$
and $h_a=1$ for $a\neq i_s, j_s$, and 
\begin{align*}
h_{i_s}&=g_{i_s}\cdots g_{i_1}x(g_{j_s}\cdots g_{j_1})^{-1}, 
\qquad s=1, \cdots, k,\\
h_{j_s}&=g_{j_s}\cdots g_{j_1}x^{-1}(g_{i_s}\cdots g_{i_1})^{-1}, 
\qquad s=1, \cdots, k.
\end{align*}
Clearly $h_{i_k}=x$ and $h_{j_k}=x^{-1}$
 by Eqn. (\ref{E:conj-pf1}).
Therefore we have the following equations for 
$s=1, 2, \cdots, k\, (mod\, k)$
\begin{equation}\label{E:conj-pf2}
h_{i_s}=g_{i_s}h_{i_{s-1}}g_{j_s}^{-1}, \quad
h_{j_s}=g_{j_s}h_{j_{s-1}}g_{i_s}^{-1},
\end{equation}
which imply that $ hs(\mu)(g)s(\rho)(h^{-1})=g$.
Note also that
$s(\rho)s(\mu)=s(\mu)s(\sigma)$ and $d(\mu)=k \, (mod\, 2)$. 
We see that the conjugation
\begin{equation}\label{E:split2}
(h, t_{\mu})(g, t_{\rho})(h, t_{\mu})^{-1}
=(hs(\mu)(g)s(\rho)(h^{-1}), z^{k}t_{\rho}^{s(\mu)})
=(g, z^{k}t_{\rho}^{s(\mu)}) ,
\end{equation}
where we used $d(\rho)=1$. 
Observe that by Eqn. (\ref{E:commut})
\begin{equation*}
t_{\rho}^{s(\mu)}=\cdots [j_1\cdots j_k][i_1\cdots i_k]\cdots
=z^{(k-1)^2}t_{\rho}=z^{k-1}t_{\rho}.
\end{equation*}
Plugging this into Eqn. (\ref{E:split2}) we obtain that
$(h, t_{\mu})(g, t_{\rho})(h, t_{\mu})^{-1}=z(g, t_{\rho})$,
and this contradiction says that each partition $\rho(c)$ must
be strict.
\end{proof}

Let $C_{\rho}$ be a conjugacy class in $\Gn$ 
of type $\rho=(\rho(c))_{c\in\G_*}\in\mathcal{P}_n(\G_*)$. 
We fix an order of the conjugacy classes of $\G$ as before:
$c^0, \cdots, c^r$. Let $T^{\rho(c^i)}$ be the special
tableau such that the numbers
$(\sum_{j=0}^{i-1}|\rho(c^j)|)+1$, $\cdots$, $\sum_{j=0}^{i}|\rho(c^j)|$ appear in the natural order
from the left to right and up to bottom in the Young
diagram of shape $\rho(c^i)$, and thus
\begin{eqnarray}\nonumber
t^{\rho(c^i)}&=&[1+a_{i-1},\cdots, \rho(c^i)_1+
a_{i-1}]\cdots \\
\label{E:specelt1}
\ && [\rho(c^i)_1+\cdots+\rho(c^i)_{l-1}+a_{i-1},\cdots, |\rho(c^i)|+a_{i-1}],
\end{eqnarray}
where $a_{i-1}=\sum_{j=0}^{i-1}|\rho(c^j)|$ and $\rho(c^i)=(\rho(c^i)_1, \cdots, \rho(c^i)_l)$.
We define the special element
$t^{\rho}$ by
\begin{equation}\label{E:specelt2}
t^{\rho}=t^{\rho(c^0)}t^{\rho(c^1)}\cdots t^{\rho(c^r)}.
\end{equation}

For each split conjugacy class $C_{\rho}$ in $\Gn$ of type 
$\rho$, we define the conjugacy class $D_{\rho}^+$ in $\tG n$
to be the conjugacy class
containing the element $(g, t^{\rho})$. We also
define $D_{\rho}^-=zD_{\rho}^+$. Then
$\theta_n^{-1}(C_{\rho})=D_{\rho}^+\cup D_{\rho}^-
$. Let $(D_{\rho}^+)^{-1}=\{x^{-1}\vert x\in
D_{\rho}^+\}$. We remark that $(D_{\rho}^+)^{-1}=z^{n(\rho')}(D_{\overline{\rho}})^+$, where
$n(\rho')$ is defined in (\ref{E:n}) and $\overline{\rho}$
is the partition-valued function given by
$\overline{\rho}(c)=\rho(c^{-1})$.

Given a partition $\lambda = (1^{m_1} 2^{m_2} \ldots )$,
we define
\[
  z_{\la } = \prod_{i\geq 1}i^{m_i}m_i!.
\]
We note that $z_{\la }$ is the order of the centralizer 
of an element of cycle-type $\la $ in $S_{|\la |}$. 

For each partition-valued function $\rho=(\rho(c))_{c\in\G_*}$
we define
$$
Z_{\rho}=\prod_{c\in\G_*}z_{\rho(c)}\zeta_c^{l(\rho(c))},
$$
which is the order of the centralizer of an element of conjugacy 
type $\rho=( \rho(c))_{ c \in \G_*}$ (see \cite{M}). 

\begin{proposition} \label{P:centralizer} The order of the
centralizer of an element of conjugacy  
type $\rho$ in $\tG n$ is given by
\begin{equation*}
\tilde{Z}_{\rho}=\begin{cases} 2Z_{\rho}, & \mbox{$C_{\rho}$ is split}\\
Z_{\rho}, & \mbox{$C_{\rho}$ is non-split.}
\end{cases}
\end{equation*}
\end{proposition}
\begin{proof} Let $C_{\rho}$ be a conjugacy class in $\Gn$.
If $\theta^{-1}(C_{\rho})$ does not split, then $\theta^{-1}(C_{\rho})$ is a conjugacy class in $\tG n$, so its
centralizer has the order $|\tG n|/|\theta^{-1}(C_{\rho})|
=|\Gn|/|C_{\rho}|=Z_{\rho}$. Otherwise 
$\theta^{-1}(C_{\rho})=D_{\rho}^+\cup D_{\rho}^-$, and
$|\tG n|/|D_{\rho}^{\pm}|=2Z_{\rho}$.
\end{proof}

\begin{theorem} \cite{HH} (1) The number of conjugacy classes of $\Gn$ and $\Gn^0$
are given by the following formulas:
\begin{align*}
|\mbox{split classes of the group $\Gn$}|&= 2 |\mathcal{SP}^1_n(\G_*)|+
|\mathcal{SP}^0_n(\G_*)|,\\
|\mbox{split classes of the group $\Gn^0$}|&= |\mathcal{SP}^1_n(\G_*)|+
2|\mathcal{SP}^0_n(\G_*)|.
\end{align*}

(2) The number of irreducible double spin
representations is equal to the number of even strict partition-valued
functions on $\G_*$, and the number of pairs of irreducible associate spin
representations is equal to the number of odd strict partition-valued
functions on $\G_*$.
\end{theorem}
\begin{proof} The first statement in Part (1) is a corollary of 
Theorem \ref{T:class} and Lemma \ref{L:euler}. To see the second equation
in (1), we observe that an 
irreducible spin representation $\pi$ decomposes
as follows when restricting to the subgroup $\tG n^0$:
\begin{equation*}
\pi|_{\tG n^0}=\begin{cases} \pi_1\oplus \pi_2, 
&\mbox{$\pi$ is a double spin,}\\
\pi, &\mbox{$\pi$ is an associate spin.}
\end{cases}
\end{equation*}
Moreover a pair of the associated spin representations,
when restricted to $\tG n^0$, become the same irreducible
representation.   
Applying the counting formulas in (1) we obtain Part (2).
\end{proof}

We remark that the number of split
conjugacy classes of $\Gn$ contained
in $\Gn^0$ is equal to $|\mathcal{SP}^1_n(\G_*)|+
|\mathcal{SP}^0_n(\G_*)|=|\mathcal{OP}_n(\G_*)|$ by Theorem \ref{T:class}.

\section{The Hopf algebra $\tRG$ of super spin characters} \label{S:hopf}
\subsection{Superalgebras and supermodules}
We basically follow the exposition of \cite{Jo} in this subsection.
A complex {\it superalgebra} 
$A=A_0\oplus A_1$ is a $\mathbb Z_2$-graded complex vector space with a binary product $A\times A\longrightarrow A$ such that $A_iA_j\subset A_{i+j}$.
 
A vector space $V=V_0\oplus V_1$ is a {\it supermodule} for a
superalgebra $A=A_0\oplus A_1$ if $A_iV_j\subset V_{i+j}$.
Elements of $V_i$ are called homogeneous.
A linear map $f: M\rightarrow N$ between two 
$A$-supermodules is a {\it super homomorphism} 
of degree $i$ if
$f(M_j)\subset M_{i+j}$ and for any homogeneous element
$a\in A$ and any homogeneous vector 
$m\in M$ we have
\[ f(am)=(-1)^{d(a)d(f)}af(m).
\]
Let 
\[
Hom_A(M, N)=Hom_A(M, N)_0\oplus Hom_A(M, N)_1,
\]
where
$Hom_A(M, N)_i$ consists of $A$-super-homomorphisms of 
degree $i$
from $M$ to $N$. 

Let $V=V_0\oplus V_1$ and $W=W_0\oplus
W_1$ be two supermodules. The tensor product $V\otimes W$ 
is also a supermodule with
$(V\otimes W)_i=\sum_{k+l=i (mod 2)}V_k\otimes W_l$.
Submodules, irreducible or simple supermodules are defined
similarly as usual. Two examples of complex simple superalgebras
are given in order.

Let $r, s\in \mathbb N$. We define $M(r|s)$ to be the $\mathbb C$-superalgebra of 
$(r+s)$-square matrices with the grading
\begin{align*}
M(r|s)_0&=\{\begin{bmatrix}
A & 0\\ 0 & D\end{bmatrix}| A\in M_{r, r}(\mathbb C), D
\in M_{s, s}(\mathbb C)\},\\
M(r|s)_1&=\{\begin{bmatrix}
0 & B\\ C & 0\end{bmatrix}| B\in M_{r, s}(\mathbb C), C
\in M_{s, r}(\mathbb C)\}\\
\end{align*}
and the operations are the underlying 
usual matrix
addition and multiplication. As in the ungraded case,
$M(r|s)$ can also be viewed as the superalgebra of 
$\mathbb Z_2$-graded linear maps
of $\mathbb C^{r|s}=\mathbb C^{r}\oplus \mathbb C^{s}$
with the usual superpositions of maps.
It is easily seen that $M(r|s)$ is a simple
$\mathbb C$-superalgebra and $\mathbb C^{r|s}$ is a simple
$M(r|s)$-supermodule.  

Another example is the $\mathbb C$-superalgebra $Q(n)$. 
As a supervector space it is defined by
\begin{align*}
Q(n)_0&=\{\begin{bmatrix}
A & 0\\ 0 & A\end{bmatrix}| A\in M_{n, n}(\mathbb C)\},\\
Q(n)_1 &=\{\begin{bmatrix}
0 & B\\ B & 0\end{bmatrix}| B\in M_{n, n}(\mathbb C)\}.
\end{align*}
The superalgebra structure is given by the usual matrix
multiplication. The space $\mathbb C^{n|n}$ is also
a $Q(n)$-supermodule under the usual matrix multiplication.

Wall \cite{Wl}
showed that these two simple superalgebras
are the only two types of simple superalgebras over $\mathbb C$. In the sequel we will call the supermodule
$\mathbb C^{r|s}$ of type $M$ if it is considered as a
$M(r|s)$-supermodule and $\mathbb C^{n|n}$ of type 
$Q$ if it is considered as a $Q(n)$-supermodule.

For any finite group $G$ and a subgroup $H$ of index $2$, 
we set the parity of elements of $H$ (resp. $G\backslash
H$) to be
even (resp. odd). The corresponding group superalgebra 
of $G$ 
is semisimple (see \cite{Jo}). In the case
of the spin wreath product $\tG n$ and the subgroup
$\tG n^0=\G^n\wr \tilde{A_n}$, this parity agrees with
the parity given by the homomorphism $d$ (see (\ref{E:parity}). As a superalgebra, 
$\mathbb C[\tG n]$ is given by 
\begin{align}
\mathbb C[\tG n]_0 &=\{\sum_g a_g g| g\in \tG n^0\},\\
\mathbb C[\tG n]_1 &=\{\sum_g a_g g| g\in \tG n^1\},
\end{align}
and the product is the usual multiplication. 

\begin{proposition} There exists an isomorphism of
$\mathbb C$-superalgebras
$$
\mathbb C[\tG n]\simeq \bigoplus_{i}M(r_i|s_i)\oplus \bigoplus
_jQ(n_j).
$$
Any finite dimensional $\mathbb C[\tG n]$-supermodule is isomorphic to a
direct sum of simple supermodules of type $M$ and $Q$.
\end{proposition}
 
By the definition of spin representations and Lemma \ref{L:euler}
we know the number of irreducible spin supermodules
for $\tG n$.
\begin{proposition} \label{P:number} 
 The number of irreducible spin supermodules of $\tG n$ is equal to 
$|\mathcal{SP}_n(\G_*)|$,
the number of strict 
partition-valued functions on $\G_*$. If $V$ is an irreducible
$\tG n$-supermodule of type $M$, then its underlying $\tG n$-module
is irreducible. If $V$ is an irreducible
$\tG n$-supermodule of type $Q$, then its underlying $\tG n$-module 
decomposes into two irreducible $\tG n$-modules $U$ and $U'$, where
$U'=U$ as a vector space and its action is given by
$a.u=(-1)^{d(a)}au$ for any homogeneous element $a\in \tG n$.
\end{proposition}

\subsection{Induced supermodules}
Let $G$ be a finite group with a central involution
$z$ and a parity epimomorphism $d: G\longrightarrow \mathbb
Z_2$ such that $d(z)=0$. Let $H$ be a subgroup of $G$ containing $z$ such that the restriction of $d$
on $H$ is not identically zero. Such a pair 
$(G, H)$ of finite groups
 will be called an 
{\it admissible pair} of finite groups.
 
The group algebras $\mathbb C[G]$ and $\mathbb C[H]$
become superalgebras with $G^0=ker(d), H^0=ker(d|_H)$
and $G^1=G\backslash G^0, H^1=G\backslash H^0$.

Let $W$ be a $\mathbb C[H]$-supermodule. We define
the {\it induced supermodule} $Ind_H^G W$ for $\mathbb C[G]$ by
\begin{equation}\label{E:ind}
Ind_H^G W=\mathbb C[G]\otimes_{\mathbb C[H]}W
\end{equation}
with the action given by $g(h\otimes w)=gh\otimes w$.
Clearly $Ind_H^G W$ is a spin supermodule
if $W$ is a spin supermodule.

The following lemma can be checked similarly as
the ordinary case \cite{Sr}.

\begin{lemma} \label{L:induction} Let $(G, H)$ 
be an admissible pair of finite groups. Let $V$ be
a 
$\mathbb C[G]$-supermodule and $W$ a $\mathbb C[H]$-sub-supermodule of $V\vert_H$. Then $V$ is equal to
$Ind_H^G W$ if and only if
\begin{equation}
V=\bigoplus_{s\in G/H}\quad sW,
\end{equation}
where $sW$ denotes the subspace $x.W$ $(x\in s)$ 
 of the supermodule $V$. 
\end{lemma}

Let $(G, K)$ be another admissible pair of finite
groups. Consider the double cosets $HsK$ of
$H$ and $K$ in $G$. For $s\in H\backslash G/K$ 
the set $H_s:=s^{-1}Hs\cap K$ is a subgroup of $K$.
 The following analog of Mackey's theorem can be proved
similarly as in the ordinary case using Lemma \ref{L:induction}.

\begin{proposition}\label{P:mackey}
Let $(G, H)$ and $(G, K)$ be admissible pairs of
finite groups as above. Then we have
\begin{equation}
Res_KInd_H^G\, W\simeq \bigoplus_{s\in H\backslash G/K}
Ind_{H_s}^K Res_{H_s} W
\end{equation}
as supermodules. 
\end{proposition}

\subsection{The space $R^-(\tG n)$}
 A {\it spin class function} on $\tG n$ 
is a class function map from $\tG n$ to
$\mathbb C$ such that
\[ f(zx)=-f(x). 
\]
Thus spin class functions vanish on non-split conjugacy classes.
A {\it spin super class function} 
on $\tG n$ 
is a spin class function $f$ on $\tG n$ such that
$f$ vanishes further on odd strict conjugacy classes. In other
words, $f$ corresponds to a complex functional on $\mathcal{OP}_n(\G_*)$ in view of
Theorem \ref{T:class}. 

Let $R^-(\tG n)$ be the $\mathbb C$-span of 
spin super class functions on $\tG n$. 
Let $R(\tG n)$ be the $\mathbb C$-span of class functions
on $\tG n$. Let $R^0(\tG n)$ be the subspace of the class functions $f(x)$ such that $f(zx)=f(x), x\in\tG n$, and let $R^1(\tG n)$ be the space
of spin class functions. Then we have
\begin{align*}
R(\tG n)&=R^0(\tG n)\oplus R^1(\tG n), \\
R^-(\tG n)&\subset R^1(\tG n), \qquad R^0(\tG n)\simeq 
R(\Gn).
\end{align*}
In this paper we will focus on the space $R^-(\tG n)$.
We remark that $R^1(\tG n)$ can be identified 
as a vector space with the Grothendieck ring
of spin representations of $\tG n$, and it is
not difficult to recover $R^1(\tG n)$ from $R^-(\tG n)$
using Proposition \ref{P:number}.

The standard inner product
$\langle \ | \ \rangle$ on $R(\tG n)$
induces an inner product on $R^-(\tG n)$. 
For two spin super class functions $f, g\in R^-(\tG n)$ we define
\begin{align} \nonumber
\langle f, g \rangle&= 
\langle f, g \rangle_{\tG n}\\ \label{E:innerprod}
&=\frac1{ | \tG n |}\sum_{x \in \tG n^0}
          f(x) g(x^{ -1})=\frac 12
\langle f, g \rangle_{\tG n^0},
\end{align}
where $\langle f, g \rangle_{\tG n^0}$ is the
inner product of $f|_{\tG n^0}$ and $g|_{\tG n^0}$
in the space of class functions on the subgroup
$\tG n^0$.
Since even split conjugacy classes of $\tG n$ have the form 
$\{D_{\rho}^+\}\cup \{D_{\rho}^-\}$ and $f(D_{\rho}^+)=-f(D_{\rho}^-)$,
we can rewrite the inner product by using Proposition \ref{P:centralizer}.
\begin{equation} \label{E:alt}
\langle f, g \rangle=\sum_{\rho\in \mathcal{OP}_n(\G_*)} 
\frac 1{Z_{\rho}}f({\rho}) 
g({\overline{\rho}}),
\end{equation}
where $f({\rho})=f(D_{\rho}^+)$ and $\overline{\rho}$ is defined in Sect. \ref{S:splitclass}.
In the sequel we will fix the value of a class function
at $\rho\in \mathcal{OP}_n(\G_*)$ to be the value at the conjugacy
class $D_{\rho}^+$. 

Let $V=V_0\oplus V_1$ be a $\tG n$-supermodule. 
We define the character of $V$ as the function 
$\chi_V: x\mapsto tr(x)$, $x\in\tG n$. 
Clearly $\chi_V(\tG n^1)=0$. 

\begin{proposition} \label{P:criterion}
The characters of irreducible spin $\tG n$-supermo\-dules form
a $\mathbb C$-basis of $R^-(\tG n)$. 
Let $\phi$ and $\ga$ be two irreducible characters
of spin supermodules, then
\begin{equation}\label{E:superorth}
\langle \phi, \ga\rangle=\begin{cases} 1 & \mbox{if $\phi\simeq \ga$, 
type $M$}\\
2 & \mbox{if $\phi\simeq \g$, type $Q$}\\
0 & \mbox{otherwise}\end{cases}.
\end{equation}
Conversely, if $\langle f, f\rangle=1$ for
$f\in R^-(\tG n)$, then $\pm f$ affords an irreducible spin $\tG n$-supermodule of type $M$. If $\langle f, f\rangle=2$, then either
$\pm f$ is the character of an 
irreducible spin supermodule of type $Q$ or
$f$ is a sum or difference of two irreducible characters of 
spin supermodules of type $M$.
\end{proposition}
\begin{proof} Let $\xi_{c}$ be the characteristic function
on the conjugacy class $c$. Then $R^-(\tG n)$ is spanned by $\xi_{c}$, where $c$ ranges over the set of split even classes. Thus $dim(R^-(\tG n))\le |\mathcal{OP}_n(\G_*)|$.

On the other hand we see that the characters of spin supermodules are class functions in $R^-(\tG n)$ since the trace of any odd endomorphism is zero. Let $\phi$ and $\ga$ be the characters of two irreducible spin supermodules of $\tG n$. It follows from Proposition \ref{P:number} that the underlying 
module of $\phi$ or $\ga$
is either irreducible module or the sum of two associated irreducible modules according to their types,
which implies immediately the orthogonality relation
(\ref{E:superorth}). Therefore the matrix of the inner product is orthogonal on the set of super spin characters.
Then by Lemma \ref{L:euler} and Proposition \ref{P:number}
$$dim(R^-(\tG n))\ge |\mathcal{SP}_n(\G_*)|=
|\mathcal{OP}_n(\G_*)|.
$$
Thus the two inequalities above become equality, and so
the  irreducible characters of spin 
$\tG n$-supermodules
form a $\mathbb Z$-basis in $R^-(\tG n)$.

The last characterization of irreducible supermodules
follows from the semi-simplicity of the superalgebra
$\mathbb C[\tG n]$ and the usual orthogonality of
ordinary irreducible characters.
\end{proof}  

\subsection{Hopf algebra structure on $\tRG$}\label{S:Hopf}
We now define one of our main objects  
\[
  \tRG = \bigoplus_{n\geq 0} R^-(\tG n). 
\]

Let $\tG n\tilde{\times}\tG m$ be the direct product
of $\tG n$ and $\tG m$ with a twisted multiplication 
\begin{equation*}
(t, t')\cdot (s, s')=(ts z^{d(t')d(s)}, t's'),
\end{equation*}
where $s, t\in\tG n, s', t'\in\tG m$ 
are homogeneous. Note that $|\tG n\tilde{\times}\tG m|=
|\tG n||\tG m|$. We define the {\it spin product}
of $\tG n$ and $\tG m$ (see \cite{HH}) by
\begin{equation}
\tG n\hat{\times}\tG m={\tG n\tilde{\times}\tG m}/
\{(1, 1), (z, z)\},
\end{equation}
which can be embedded into the spin group $\tG{n+m}$
canonically by letting 
\begin{equation}\label{E:embed}
((g, t_i'), 1)\mapsto (g, t_i), \qquad (1, (g, t_j'')) \mapsto (g, t_{n+j}),  
\end{equation}
where $i=1, \ldots, n-1$, $j=1,\ldots, m-1$. We will
identify $\tG n\hat{\times}\tG m$ with its image in
$\tG {m+n}$ and regard it as a subgroup of $\tG{m+n}$. 
Clearly $\theta_{n+m}(\tG n\hat{\times}\tG m)$ is
the pull-back of $\Gn\times\Gm$.

\begin{remark} \label{R:embed}
 Partition $\{1, 2, \cdots, n+m\}$
into a disjoint union of subsets $I$ and $J$ with
$|I|=n$ and $|J|=m$. Then $\tG n\hat{\times}\tG m$ can be
embedded into $\tG {n+m}$ using a similar map as (\ref{E:embed}) by mapping 
the generators of $\tS n$ and $\tS m$ to
the generators of $\tS {n+m}$ indexed by $I$ and $J$ 
respectively. One can check that all such embeddings
of $\tG n\hat{\times}\tG m$ are conjugate subgroups
in $\tG{n+m}$.
\end{remark}  

The subgroup $\tG n\hat{\times}\tG m$ has a distinguished subgroup
of index $2$ consisting of even elements given by $d$. We
define $R^-(\tG n\hat{\times}\tG m)$ to be the space
of spin class functions on $\tG n\hat{\times}\tG m$
that vanish on odd conjugacy classes of $\tG n\hat{\times}\tG m$.

For two spin supermodules $U$ and $V$ of
$\tG n$ and $\tG m$ we define the {\it super (outer)-tensor product} $U\otimes V$
by
\begin{equation*}
(t, s)\cdot (u\otimes v)=(-1)^{d(s)d(u)}(tu\otimes sv),
\end{equation*}
where $s$ and $u$ are homogeneous elements. We see
immediately that
\begin{align*}
(z', z'')\cdot (u\otimes v)&=(-u)\otimes (-v)=u\otimes v,\\
(z, 1)\cdot (u\otimes v)&=-(u\otimes v).
\end{align*}
This says that $U\otimes V$ is a spin 
$\tG n\hat{\times}\tG m$-supermodule.

The following is a direct generalization of a result in \cite{Jo}
for trivial $\G$.

\begin{proposition}\label{P:supermult}
Let $U$ and $V$ be simple 
supermodules for $\tG n$ and $\tG m$ respectively. Then
\item 1) If both $U$ and $V$ are of type M, then $U\otimes V$
is a simple $\tG n\hat\times\tG m$-supermodule of type M.
\item 2) If $U$ and $V$ are of different type, then
$U\otimes V$
is a simple $\tG n\hat\times\tG m$-supermodule of type Q.
\item 3) If both $U$ and $V$ are type Q, then $U\otimes V
\simeq N\oplus N$ for some simple $\tG n\hat\times\tG m$-supermodule 
$N$ of type M.
\end{proposition}

The super (outer)tensor product defines
an isometric isomorphism 
\begin{equation}\label{E:twprod}
R^-(\tG n ) \bigotimes R^-(\tG m) 
 \stackrel{\phi_{n, m}}\rightarrow R^-(\tG n \hat{\times} \tG m),
\end{equation}
which is actually an isomorphism over $\mathbb Q$ by Proposition \ref{P:supermult}.

The space $\tRG$ carries a multiplication defined by the composition 
\begin{equation}\label{E:hopf1}
 m: R^-(\tG n ) \bigotimes R^-(\tG m) 
 \stackrel{\phi_{n, m} }{\longrightarrow} R^-(\tG n \hat{\times} \tG m)
 \stackrel{Ind}{\longrightarrow} R^-( \tG{n + m}),
\end{equation}
and a comultiplication defined by the composition
\begin{equation}\label{E:hopf2}
\begin{aligned}
\Delta: R^-(\tG n ) &\stackrel{Res}{\longrightarrow}
 \bigoplus_{m =0}^n R^-( \tG {n - m} \hat{\times} \tG m)\\
 &\stackrel{\phi^{-1}}{\longrightarrow}
 \bigoplus_{m =0}^n R^-( \tG {n - m}) \bigotimes 
R^-(\tG m).
\end{aligned}
\end{equation}
Here $Ind$ and $Res$ denote the induction (see 
(\ref{E:ind}) and restriction functors respectively.
The isomorphism $\phi^{-1}$ is equal to 
$\oplus_{0\leq m\leq n}\phi^{-1}_{n-m, m}$ (see (\ref{E:twprod})).

\begin{theorem}
The above operations define a Hopf algebra structure for $\tRG$.
\end{theorem} 
\begin{proof} Using Remark \ref{R:embed} twice we observe that the following two embeddings
give rise to two conjugate subgroups in $\tG{n+m+l}$
(see Remark \ref{R:embed}):
\begin{equation*}
(\tG n\hat{\times}\tG m)\hat{\times}\tG l\hookrightarrow
\tG{n+m+l}\hookleftarrow \tG n\hat{\times}
(\tG m\hat{\times}\tG l).
\end{equation*}
Using this and Lemma \ref{L:induction} we can easily check the associativity of the product.

For a simple supermodule $V$ we define
\begin{equation}  
c(V)=\begin{cases} 0 & \mbox{$V$ is type M}\\
  1 & \mbox{$V$ is type Q}\end{cases}.
\end{equation} 
Let $V_1$, $V_2$, and $V_3$ be simple supermodules
for $\mathbb C[\tG n]$, $\mathbb C[\tG m]$, and $\mathbb C[\tG l]$ respectively. It is easy to see that
$c(V_1, V_2)=c(V_1)c(V_2)$ satisfies the cocycle condition
\begin{equation} \label{E:cocycle}
c(V_1, V_2)+c(V_1\otimes V_2, V_3)
=c(V_2, V_3)+ c(V_1, V_2\otimes V_3).
\end{equation}
Therefore we can define $c(V_1\otimes V_2\otimes V_3)$ 
to be either
of the above expressions. 

Using the cocycle $c$ we prove the coassociativity as follows.
Let $U$ be a $\mathbb C[\tG n]$-supermodule and 
suppose that 
$Res_{\tG m\hat{\times}\tG l\hat{\times}\tG k}U=\oplus_{i}U_i(m, l, k)$ as 
an irreducible decomposition. Then we have
\begin{align*}
(1\otimes \Delta)\Delta(U)&=(1\otimes\phi^{-1})\phi^{-1}
\bigoplus_{m+l+k=n} Res_{\tG m\hat{\times}\tG l\hat{\times}\tG k}U\\
&=\bigoplus_{m+l+k=n, i} 2^{-c(U_i(m, l, k))}U_i(m, l, k)\\
&=(\phi^{-1}\otimes 1)\phi^{-1}
\bigoplus_{m+l+k=n} Res_{\tG m\hat{\times}\tG l\hat{\times}\tG k}U\\
&=(\Delta\otimes 1)\Delta(U),
\end{align*}
where we used the cocycle condition in the third equation, and the notation
$\sum 2^{c(U_i)}U_i$ stands for the multiplicity-free summation of the irreducible components (c.f. Proposition
\ref{P:supermult} and definition of $\phi$).

Finally we look at the compatibility of multiplication and comultiplication. Fix $m$ and $n$,
it follows from Proposition \ref{P:mackey} that
\begin{align*}
\Delta(U\cdot V)&=\bigoplus_{k+l=m+n}\phi^{-1}_{k, l}
Res_{\tG k\hat{\times}\tG l}Ind_{\tG m\hat{\times}\tG n}^{\tG{m+n}}
\phi_{m, n}(U\otimes V)\\
&=\bigoplus_{k+l=m+n}\bigoplus_{s}\phi^{-1}_{k, l}
Ind_{(\tG m\hat{\times}\tG n)_s}^{\tG k
\hat{\times}\tG l}
Res_{(\tG m\hat{\times}\tG n)_s}\phi_{m, n}(U\otimes V)_s,
\end{align*}
where $s$ runs through the double cosets 
$\tG m\hat{\times}\tG n\backslash
\tG{m+n}/\tG k\hat{\times}\tG l$. Notice that the
double cosets are in one-to-one correspondence with the double cosets $\Gm{\times}\Gn\backslash
\G_{m+n}/\G_k{\times}\G_l$. Again by the cocycle property 
of $c$ and counting the double cosets we can check that the last
summation is 
exactly $\Delta(U)\cdot\Delta(V)$.
\end{proof}

\begin{remark} Our Hopf algebra 
is different from
that of \cite{HH} where a bigger space
than our $\tRG$ was used. 
\end{remark}

The standard bilinear form
in $\tRG$ is defined in terms of those on $R^-( \tG n )$
as follows:
\[
\langle u, v \rangle
 = \sum_{ n \geq 0} \langle u_n, v_n \rangle_{\tG n},
\]
where
$u = \sum_n u_n$ and $v = \sum_n v_n$ with $u_n, v_n\in \tG n$.
\section{Basic spin representations of $\tG n$}
%
%
%
%
%
%
\subsection{A weighted bilinear form on $R(\G)$ and $R^-(\tG n)$}
In \cite{FJW1} we introduced the notion of weighted bilinear forms
on $R(\G)$ and coherently combined several examples in this concept.
We will also similarly define weighted bilinear forms on
the space $R^-(\tG n)$.

Let $\wt$ be a self-dual class function in $R(\G)$, i.e.
$\wt(c)=\wt(c^{-1})$. Let $*$ denote the product of two
characters in $R(\G)$, which is afforded by the tensor product. 
Let $a_{ij} \in \mathbb C$ be the (virtual) multiplicities 
of $\g_j$ in $ \wt * \g_i $: 
\begin{eqnarray}  \label{eq_tens}
 \wt  *  \g_i 
  = \sum_{j =0}^r a_{ij} \g_j.
\end{eqnarray}
We denote further by $A$ the $ (r +1) \times (r +1)$ matrix 
$ ( a_{ij})_{0 \leq i,j \leq r}$. Then the weighted
bilinear form $\langle f, g \rangle_{\wt }$ is defined by
$$
  \langle f, g \rangle_{\wt } = \langle \wt * f ,  g \rangle_{\G },
   \quad f, g \in R( \G).
$$ 
Alternatively it can be explicitly given by
\begin{eqnarray}
  \langle f, g \rangle_{\wt } 
   & =& \frac 1{ |\G|} \sum_{ x \in \G} \wt (x)f(x) g (x^{ -1}) 
     \nonumber  \\ \label{eq_twist1}
   & =& \sum_{c \in \G_*} \zeta_c^{ -1} \wt (c) f(c) g (c^{ -1})\\
& =& \sum_{c \in \G_*} \zeta_c^{ -1} \wt (c) f(c^{-1}) g (c)
     \label{eq_twist2}.
\end{eqnarray}
In particular, Eqn.~(\ref{eq_tens}) is equivalent to
\begin{equation}\label{E:matrixA}
   \langle \g_i, \g_j \rangle_{\wt } = a_{ij}.
\end{equation}
The self-duality implies that $A$ is a symmetric matrix. 
Note that the weighted bilinear
form becomes the standard bilinear form when $\xi=\g_0$, the trivial character of $\G$.

Let $V$ be a spin supermodule for $\tG n$ and $W$ a module
for $\Gn$.  As a $\mathbb Z_2$-graded
vector space 
$W\otimes V=W\otimes V_0\oplus W\otimes V_1$ and the action of $\tG n$
is defined by
\begin{equation}\label{E:circtensor}
(g, z^pt_{\rho})(w\otimes v)=(g, s(\rho))
\cdot w\otimes (g, z^pt_{\rho})\cdot v, \qquad g\in \G^n, \sigma\in\mathcal{P}_n(\G_*).
\end{equation}
It is easy to check that the tensor product
$V\otimes W$ is a spin
$\tG n$-supermodule.
This construction defines
a morphism: 
\begin{equation}\label{E:tensorprod}
R(\Gn)\otimes R^-(\tG n)\stackrel{*}{\longrightarrow} R^-(\tG n).
\end{equation}

Let us recall the construction of character $\eta_n(\xi)$ in \cite{W, FJW1}.
Let $\gamma$ be an irreducible character of $\G$ afforded by
the $\G$-module $V$, the tensor product $V^{\otimes n}$ is naturally
a $\Gn$-module by the direct product action of $\G^n$
composed with permutation action of 
the symmetric group $S_n$. The resulting character of $\Gn$ is
denoted by $\eta_n(\ga)$. Furthermore we can extend $\eta_n$
from $\G^*$ to $R(\G)$. The character value of $\eta_n(\xi)$
at the class $\rho=(\rho(c))$ is given by
\begin{equation}\label{E:weightchar}
\eta_n(\xi)(\rho)=\prod_{c\in\G_*}\xi(c)^{l(\rho(c))}.
\end{equation}
It is clear that the class function $\eta_n(\xi)$ is self-dual
as long as $\xi$ is.

We now introduce a {\em weighted bilinear form} on $R^-( \tG n)$ by letting
$$
  \langle  f, g\rangle_{\wt, \tG n } =
   \langle \eta_n (\wt )* f, g \rangle_{\tG n} ,
   \quad f, g \in R^-( \tG n),
$$
where we used the map (\ref{E:tensorprod}).
The self-duality of
$\eta_n (\wt)$ implies that the 
bilinear form $\langle \ , \ \rangle_{\wt}$ is symmetric.

\begin{remark}
 When $n =1$, this weighted bilinear form obviously reduces
 to the weighted bilinear form on $R( \G)$ defined in
(\ref{eq_twist1}-\ref{eq_twist2}). 
\end{remark}

The bilinear form on $\tRG = \bigoplus_{n} R^-(\tG n)$ is given by
\[
\langle u, v \rangle_{\wt}
 = \sum_{ n \geq 0} \langle u_n, v_n \rangle_{\wt, \tG n } ,
\]
where
$u = \sum_n u_n$ and $v = \sum_n v_n$ with $u_n, v_n\in R^-(\tG n)$.
\subsection{Basic spin representations}\label{S:basicspin}
Let the Pauli spin matrices be
\begin{align*}
\sigma_0&=\begin{bmatrix} 1 & 0\\ 0 & 1\end{bmatrix},
\qquad \sigma_1=\begin{bmatrix} 0 & 1\\ 1 & 0\end{bmatrix},\\
\sigma_2&=\begin{bmatrix} 0 & -\sqrt{-1}\\ \sqrt{-1}  & 0\end{bmatrix},
\qquad \sigma_3=\begin{bmatrix} 1 & 0\\ 0& -1\end{bmatrix}.
\end{align*}
Let $C_{2k}$ be the Clifford algebra generated by
$e_1, e_2, \cdots e_{2k}$ with relations:
\begin{equation}\label{E:clifford}
\{e_i, e_j\}=e_ie_j+e_je_i=-2\delta_{ij}.
\end{equation}
Thus $e_j^2=-1$. The Clifford algebra $C_{2k}$ is endowed with a natural superalgebra
structure by letting the parity of $e_i$ to be odd for each
$i$.
When $k=1$, one has that $C_2\simeq M(1|1)=End(\mathbb C^{1|1})$
and the action of $C_2$ on $\mathbb C^{1|1}$ is given 
by the Pauli 
spin matrices: 
\[ 
e_1\mapsto \sqrt{-1}\sigma_1, \qquad 
e_2\mapsto \sqrt{-1}\sigma_2.
\]

More generally we have $C_{2k}=End(\otimes^k\mathbb C^{1|1})\simeq M(2^{k-1}|2^{k-1})$.
The tensor product $\otimes^k \mathbb C^{1|1}$ admits
a canonical supermodule structure for the Clifford algebra $C_{2k}$
under the action
\begin{align} \label{E:spinaction1}
e_{2j-1}&\longrightarrow \sqrt{-1}\sigma_3^{\otimes(j-1)}
\otimes \sigma_1\otimes \sigma_0^{\otimes(k-j)}, \quad j=1, \ldots, k,\\ \label{E:spinaction2}
e_{2j}&\longrightarrow \sqrt{-1}\sigma_3^{\otimes(j-1)}
\otimes \sigma_2\otimes \sigma_0^{\otimes(k-j)},
\quad j=1, \ldots, k.
\end{align}
The above formulas define explicitly the structure
of a simple $C_{2k}$-supermodule on $\otimes^k \mathbb C^{1|1}$.

Let $C_{2k+1}$ be the Clifford algebra generated by
$e_i$, $i=1, \ldots, 2k+1$ with similar relations like
(\ref{E:clifford}). We embed $C_1$ into $End(\mathbb C^{1|1})$ by $1\mapsto Id$, $e_1\mapsto \sqrt{-1}\sigma_1$.
Then
\begin{equation*}
C_{2k+1}\simeq C_{2k}\otimes C_1 \hookrightarrow 
End(\otimes^{k+1}\mathbb C^{1|1})
\end{equation*}
gives a $C_{2k+1}$-supermodule structure on  
$\otimes^{k+1}\mathbb C^{1|1}$. The explicit
action is given by the same formulae 
(\ref{E:spinaction1}-\ref{E:spinaction2}), except
that $j=1, \cdots, k+1$ in (\ref{E:spinaction1}).
Observe that $C_{2k+1}\simeq Q(2^k)$.

It is well-known (see e.g. \cite{Jo}) that there exists
an embedding of $\tS n$ into the multiplicative Clifford group
of units in $C_{n-1}$. Therefore 
$\otimes ^{[\frac n2]}\mathbb C^{1|1}$ can be regarded as a
$\tS n$-supermodule, which is called
the {\it basic spin supermodule} for $\tS n$.
More explicitly we have 

\begin{proposition}\cite{Sc, Jo} The basic spin
supermodule 
for $\tS n$ is $\otimes^{[\frac n2]}\mathbb C^{1|1}$ with the action
\begin{equation}\label{E:spinaction}
t_j\mapsto \sqrt{\frac{j+1}{2j}}e_j-\sqrt{\frac{j-1}{2j}}e_{j-1}, \qquad j=1, \cdots, n-1.
\end{equation}
Here we take $e_0=0$. Its character $\chi_n$ is given by
\begin{equation}\label{E:spinchar}
\chi_n(\alpha)=\begin{cases} 2^{l(\alpha)/2} &\mbox{if
$\alpha\in \mathcal{OP}_n$, $n$ even}\\
2^{(l(\alpha)-1)/2} &\mbox{if
$\alpha\in \mathcal{OP}_n$, $n$ odd}\\
0 & \mbox{otherwise.}
\end{cases}
\end{equation}
In particular $\chi_n(1)=2^{\lceil\frac n2\rceil}$. Here
$\lceil a\rceil$ denotes the largest integer $\leq a$.
\end{proposition}

\begin{proposition} \cite{Sc, Mo1} \label{P:spinaction2}
1) Let $n\geq 1$ be an odd 
integer. The basic spin supermodule $\otimes^{(n-1)/2}\mathbb C^{1|1}$
is an irreducible $\tS n$-module under the action
(\ref{E:spinaction}). Its character $\chi_n$ is given
by the second equation of (\ref{E:spinchar}). In particular 
$\chi_n(1)=2^{(n-1)/2}$.

2) Let $n\geq 1$ be an even integer. The basic
spin supermodule is a reducible $\tS n$-module under the action
(\ref{E:spinaction}) and decomposes into two
irreducible $\tS n$-modules whose characters $\chi_n^{\pm}$ 
are given by 
\begin{equation}\label{E:spinchar1}
\chi_n^{\pm}(\alpha)=\begin{cases} 2^{(l(\alpha)-2)/2} 
&\mbox{if $\alpha\in \mathcal{OP}_n$,}\\
\pm{(\sqrt{-1})^{n/2}}\sqrt{\frac n2} & \mbox{if
$\alpha=(n)$,}\\
0 & \mbox{otherwise.}
\end{cases}
\end{equation}
In particular, $\chi_n^{\pm}(1)=2^{(n-2)/2}$.
\end{proposition}
\subsection{The spin character $\pi_n(\g)$ of $\tG n$ 
\label{S:spinrep}}
Let $V$ be a $\G$-module afforded by the character
$\g \in R(\G)$, and let $U$ be a spin supermodule (resp. module)
of $\tS n$ with the character $\pi$. The tensor product 
$V^{ \otimes  n}\otimes U$ has a canonical
spin supermodule (resp. module) structure for $\tG n$ 
as follows 
(compare (\ref{E:circtensor})).
For any $g=(g_1, \cdots, g_n)\in \G^n$ let
$(g, z^pt_{\rho})$ be an element in $\tG n$.
The supermodule (resp. module) structure
is defined by
\begin{align*}
(g, z^pt_{\rho}).&(v_i\otimes\cdots\otimes v_n\otimes u)\\
&=g_1v_{s(\rho)^{-1}(1)}\otimes \cdots g_nv_{s(\rho)^{-1}(n)}
\otimes (z^pt_{\rho}u).
\end{align*}
We denote by $\pi_n(\g)$ the character
of the constructed spin supermodule (resp. module). 

Recall that the conjugacy class $D^+_{\rho}$ contains an element
$(g, t^{\rho})$ (see (\ref{S:conjclass})).

\begin{proposition}\label{P:superchar1} Let $\pi$ be the character of a spin
$\tS n$-supermodule. Then the character values of ${\pi_n}(\g)$
at the conjugacy classes $D^{\pm}_{\rho}$ ($\rho\in \mathcal{OP}_n(\G_*)$)
are given by
\begin{equation} \label{eq_term}
\pi_n(\g)(D_{\rho}^{\pm})=\pm\pi(t^{\rho})\prod_{c\in\G_*}
\g(c)^{l(\rho(c))}.
\end{equation}
\end{proposition}
\begin{proof} 
 Consider $(g, z^pt_{\rho}) \in \tG n,$ where 
 $g =(g_1, \ldots, g_n) \in \G^n$
 and $t_{\rho}$ is an $n$-cycle, say $t_{\rho} = [12 \ldots n]$. 
 Denote by $e_1, \ldots, e_k$ a basis of $V_{\g}$, and we write
 $ g e_j = \sum_i a_{ij} (g) e_i $,
 $ a_{ij} (g) \in \mathbb C .$
 It follows that
 \begin{align*}
  (g, z^pt_{\rho})& (e_{j_1} \otimes \ldots \otimes e_{j_n}\otimes u)\\
   &= g_1 (e_{j_n}) \otimes g_2(e_{j_1}) 
       \ldots \otimes g_n (e_{j_{n -1}})\otimes z^pt_{\rho}(u),
 \end{align*}
 
 Thus we obtain
 \begin{eqnarray*}
  \pi_n (\g) (z^pt_{\rho})
  & =& \mbox{trace } a(g_n) a(g_{n-1}) \ldots a(g_1)\pi(z^pt_{\rho})   \\
  & =& \mbox{trace } a(g_n g_{n -1} \ldots g_1)\pi(z^pt_{\rho}) = \g (c)
\pi(z^pt_{\rho}),
 \end{eqnarray*}
 where we notice that 
 $g_n g_{n -1} \ldots g_1 $ lies in $c \in \G_*$.
 
 Given $x \hat{\times} y \in \tG n$, where
 $x \in \tG r$ and $y \in \tG{n -r}$, we clearly have 
 $ \pi_n (\g ) (x \hat{\times} y) = \pi_n (\g) (x) \pi_n (\g)(y).$
 Thus it follows that for the conjugacy class $D_{\rho}^+\in\tG n$ 
of type $\rho$, we have
\begin{equation*} 
\pi_n(\g)(D^{\pm}_{\rho})=\pm\pi(t^{\rho})\prod_{c\in\G_*}
\g(c)^{l(\rho(c)}, 
\end{equation*}
 where $|| \rho || =n$. 
\end{proof}

Since the sign character is trivial at even classes, 
we can extend naturally $\pi_n$ to 
a map from $R(\G)$ to $R^-(\tG n)$ (compare with \cite{W, FJW1}).
For two $\G$-characters $\beta$ and $\g $ we define.
\begin{eqnarray}  \label{eq_virt}
  \pi_n (\beta - \g) =
  \sum_{m =0}^n ( -1)^m Ind_{\tG {n -m} \hat{\times} \tG m }^{\tG n}
   [ \pi_{n -m} (\beta) \otimes \pi_m (\g ) ].
\end{eqnarray}

When $n$ is even, the character  
$\chi_n$ of the basic spin supermodule
(see Sect. \ref{S:basicspin}) decomposes into the sum of irreducible characters $\chi^{\pm}_n$ of $\tG n$-modules. 
For each $c\in\G_*$, we define the special partition-valued function $c^{(n)}\in\mathcal{P}(\G_*)$ such that
\begin{equation} \label{E:speclass}
c^{(n)}(c)=(n), \qquad c^{(n)}(c')=\emptyset,\quad
\mbox{for
$c'\neq c$}.
\end{equation}
The following corollary is an immediate consequence of Propositions \ref{P:superchar1} and \ref{P:spinaction2}.
\begin{corollary} \label{C:char}
(1) The character
value of
$\chi_n(\g)$ at the conjugacy class $D_{\rho}^+$ of type $\rho$ is 
\begin{equation}\label{E:charvalue}
\chi_n(\g)(\rho)=\begin{cases}2^{(l(\rho)-\overline n)/2}
\prod_{c\in\G_*}
\g(c)^{l(\rho(c))} & \rho\in \mathcal{OP}_n(\G_*)\\
0 & \mbox{otherwise}\end{cases},
\end{equation}
where $\overline n$ is $0$ or $1$ depending on whether
$n$ is even or odd.

(2) Let $n$ be an even positive integer.
The character
values of 
$\chi_{n}^{\pm}(\g)$ at the conjugacy class
$D_{\rho}^+$ of type $\rho$ are 
\begin{equation}
\chi_n^{\pm}(\g)({\rho})=\begin{cases}2^{(l(\rho)-2)/2}
\prod_{c\in\G_*}
\g(c)^{l(\rho(c))} & \rho\in \mathcal{OP}_n(\G_*) \\
\pm(\sqrt{-1})^{n/2}\sqrt{\frac n2}
\g(c) & \rho=c^{(n)} \\
0 & \mbox{otherwise}\end{cases}.
\end{equation}
\end{corollary}

\subsection{Two specializations}  \label{subsec_mcka}
Let $d_i = \g_i (c^0)$ be the dimension of 
the irreducible representation of $\G$ afforded
by the character $\g_i$. Let $A$ be the matrix of
the bilinear form $\langle \ \ | \ \ \rangle$
on $R(\G)$
with respect to the basis $\g_i$.
Observe  that  
the vector
 $$
   v_i = ( \g_0 (c^i), \g_1 (c^i), \ldots, \g_r (c^i ) )^t ,
\qquad ( i =0 , \ldots, r)
 $$
 is an eigenvector of the matrix $A$
 with eigenvalue $ \wt (c^i )$. 
 
Two special choices of the weight function $\xi$ will be
our prototypical examples. The first choice is that $\xi=\g_0$, the trivial character.

Let $\pi$ be the character of
the $2$-dimensional representation of $\G$
given by the embedding of $\G$ in $SL_2(\mathbb C)$.
Let
$$
  \wt = 2 \g_0 - \pi.
$$
Then the weighted bilinear form 
$\langle \ , \ \rangle_{\wt}$ on $\RG$
becomes positive semi-definite. The radical of this
bilinear form is one-dimensional and span\-ned by 
the character of the regular representation of $\G$
$$ 
  \delta = \sum_{i = 0}^r d_i \g_i.
$$
The following is the well-known list of finite subgroups of 
$SL_2(\mathbb C)$:
the cyclic, binary dihedral, tetrahedral, octahedral 
and icosahedral groups. Mc\-Kay observed that they are in one-to-one correspondence to simply-laced Dynkin diagrams
of affine types \cite{Mc}: 
$a_{ii} =2$ for all $i$;
if $ \G \neq \mathbb Z / 2 \mathbb Z$ 
and $ i \neq j$ then $a_{ij} = 0 $ or $-1$. If 
$ \G=\mathbb Z / 2 \mathbb Z$ then $a_{01} = -2$.

\section{Twisted Heisenberg algebras and $\tG n$} \label{sect_heis}
\subsection{Twisted Heisenberg algebra $\thg $}
      Let $\thg $ be the infinite dimensional
Heisenberg algebra over $\mathbb C$, associated with 
a finite group $\Gamma$ and a self-dual class function $\wt\in R(\G)$,
with generators $a_m(\gamma), m \in 
2\Z+1, \gamma \in\G^*$ 
and a central element $C$ subject to the relations:
\begin{equation}  \label{eq_heis}
[a_m( \gamma), a_n(\gamma ')]
 = \frac m2 \delta_{m, -n} \langle \gamma, \gamma' \rangle_{\wt } C, 
 \quad m, n \in 2\mathbb Z+1, \, \g, \g ' \in \G^*.
\end{equation}
We extend $a_m (\g )$ to all
$\g = \sum_{ i =0}^r s_i \g_i \in R(\G )$ $(s_i \in \mathbb C)$
by linearity:
$ a_m ( \g ) = \sum_i s_i \, a_m (\g_i )$. 

The Heisenberg algebra may contain a large center
because the bilinear form $\langle \cdot , \cdot \rangle_{\wt }$ 
may be degenerate. 
The center of $\thg $ is spanned by $C$ together with
$a_m ( \g),  m \in 2\mathbb Z+1, \g \in R^0$,  
the radical of the bilinear form 
$\langle \cdot , \cdot \rangle_{\wt }$ in $R(\G)$. 

For $m\in 2\mathbb Z+1, c \in \G_*$ we introduce another basis for
$\thg$:
\begin{equation}\label{E:classba}
 a_{ m}( c) = \sum_{ \g\in \G^*} \gamma(c^{-1}) a_m( \g ).
\end{equation}

 The orthogonality of the irreducible
characters of $\G $ (\ref{eq_orth}) implies that
\begin{eqnarray*}
  a_m( \g ) 
   = \sum_{c \in \G_*} \zeta_c^{ -1}
   \g  (c) a_m(c).
\end{eqnarray*}

\begin{proposition}  \label{prop_orth}
 The commutation relations among the new basis for $\hg$
 are given by
 \begin{eqnarray*}
  [ a_m( {c'}^{ -1}), a_n( c )]
    & =& \frac m2 \delta_{m, -n}\delta_{c', c} \zeta_c \wt (c) C,
   \quad c, c' \in \G_*, 
\end{eqnarray*}
where $m, n\in 2\mathbb Z+1$.
\end{proposition}

\begin{proof} The proof is similar to the 
untwisted case \cite{FJW1}.
\end{proof}
\subsection{Action of $\thg$ on $\tSG$ and $\tSGG$}
     Denote by $\tSG $ the symmetric algebra 
generated by $a_{-n}(\g), n \in 2\mathbb Z_++1, \g\in \Gamma^*$.
There is a natural degree operator on $\tSG $
$$
  \deg (a_{ -n}( \g)) = n , \qquad n\in 2\mathbb Z_++1, 
$$ 
which makes $\tSG $ into a $\mathbb Z_+$-graded space. 

We define an action of $\thg$ on $\tSG$ as follows:
$a_{-n}( \g), n >0$ acts as     
a multiplication operator on $ \tSG $ and $C$ as the identity
operator; $a_n (\g),$ $ n>0$ acts as a derivation of the symmetric
algebra
\begin{eqnarray*}
  a_n (\g). a_{-n_1}( \alpha_1) a_{-n_2} (\alpha_2)
    \ldots a_{-n_k}( \alpha_k)   \\
 = \sum_{i =1}^k \delta_{n,n_i}
 \langle \g , \alpha_i \rangle_{\wt }
   a_{-n_1}( \alpha_1) a_{-n_2}(\alpha_2) \ldots 
  \check{a}_{-n_i}( \alpha_i) \ldots a_{-n_k}(\alpha_k )  .
\end{eqnarray*}
Here $n_i > 0, \alpha_i \in R(\G)$ for $i =1, \ldots , k$,
and $\check{a}_{-n_i}( \alpha_i)$ means the very term 
is deleted. In other word, the operator
$a_n (\g ), n > 0, \g \in R^0$ acts as $0$, and
$a_n (\g ), n > 0, \g \in R(\G) - R^0$ acts as certain non-zero
differential operator. Note that 
$\tSG$ is not an irreducible
representation over $\hg$ in general since the bilinear
form $\langle \ , \ \rangle_{\wt}$ may be degenerate.

Denote by $\SGO$ the ideal in the symmetric algebra $ \tSG$
generated by $a_{-n}(\g), n \in \mathbb N, \g \in R^0$. 
Denote by $\tSGG$ the quotient $\tSG /\SGO$.
It follows from the definition that $\SGO$ is a subrepresentation
of $\tSG$ over the Heisenberg algebra $\thg$. 
In particular, this induces a Heisenberg algebra action on $\tSGG$
which is irreducible. The unit $1$ in the
symmetric algebra $\tSG$ is the highest weight vector. 
We will also
denote by $1$ its image in the quotient $\tSGG$. 
\subsection{The bilinear form on $\tSG $}
       The space $ \tSG $ admits a bilinear form
$\langle \ ,  \ \rangle_{\wt } '$ determined by
\begin{eqnarray}  \label{eq_bili}
 \langle 1, 1 \rangle_{\wt}' = 1, \quad
 a_n(\g)^* = a_{-n}(\g), \qquad n\in 2\mathbb Z+1.
\end{eqnarray}
Here $a_n(\g)^*$ denotes the adjoint of $a_n(\g)$.

  For any partition $\la =( \la_1, \la_2, \dots)\in
\mathcal{OP}$ and $\g \in \G^*$, 
we define 
$$
  a_{-\la}( \g) = a_{-\la_1}( \g)a_{ - \la_2}( \g) \dots .
$$
For $\rho = ( \rho (\g) )_{ \g \in \G^*} \in {\mathcal
{OP}}(\G^* )$, 
we define 
$$
  a_{ - \rho} = \prod_{\g \in \G^*}  a_{ - \rho(\g)}(\g).
$$
It is clear that $a_{ - \rho},
\rho \in {\mathcal {OP}}(\G^* )$ form a basis for $\tSG $. 

Similarly we define 
\begin{eqnarray*}
  a_{ - \la} (c )  = a_{ - \la_1}(c) a_{ - \la_2} (c) \ldots
\end{eqnarray*}
for any partition $ \la = ( \la_1, \la_2, \ldots )\in
\mathcal{OP}$
and $c \in \G_*$.
For any $\rho = ( \rho (c) )_{ c \in \G_* } \in 
 \mathcal {OP} ( \G_* )$, we further define
\begin{eqnarray*}
  a_{- \rho}' = \prod_{ c \in \G_*} a_{ - \rho (c)} (c).
\end{eqnarray*}
The elements $a_{ - \rho} ', 
\rho \in {\mathcal {OP}}(\G_* )$ provide a new
$\mathbb C$-basis for $\tSG $.

Recall that
$\overline{\rho}  \in  \mathcal {OP} ( \G_* )$ is given by
assigning to $c \in \G_*$ the partition $\rho (c^{-1})$,
which is the composition of $\rho$ with the involution
on $\G_*$ given by $c \mapsto c^{-1}$. It follows from
Proposition~\ref{prop_orth} that

\begin{eqnarray}  \label{eq_inner}
  \langle a_{ - \rho'}', a_{ - \overline{\rho} }' \rangle_{\wt }'
  = \delta_{\rho ', \rho }
    \frac{Z_{\rho}}{2^{l(\rho)}} 
\prod_{c \in \G_*} \wt (c)^{l (\rho (c))},
 \quad \rho ', \rho  \in \mathcal{OP}(\G_*).
\end{eqnarray}

\begin{remark}
  $\SGO $ can be characterized as the radical of the
 bilinear form $\langle \ , \ \rangle_{\wt}'$ in $\tSG$. 
 Thus the bilinear form $\langle \ , \ \rangle_{\wt} '$
 induces a bilinear form on $\tSGG$ which will be denoted
by the same notation. 
\end{remark}
\section{Isometry between $\tRG$ and $\tSG$} 
\label{sect_isom}
\subsection{The characteristic map $\ch$}
We define a $\mathbb C$-linear map 
$ch: \tRG \longrightarrow \tSG$ by letting
\begin{equation}\label{E:ch}
ch (f) 
= \sum_{\rho \in \mathcal {OP}(\G_*)}\frac{2^{l(\rho)/2}}{Z_{\rho}}
f(\rho) 
a_{-\overline{\rho}}',
\end{equation}
where $f(\rho)=f(D_{\rho}^+)$.
The map $ch $ is called the {\em characteristic map} (compare
with \cite{Sc, Jo} for $\G$ trivial). 

Fix $n\in 2\mathbb Z_++1$ in this paragraph. Denote by $D_n(c)^{+},  (c \in \G_*)$ the conjugacy class in $\tG n$
of elements $(x, t_s) \in \Gn$ such that $s$ is an 
$n$-cycle and the cycle product of $(x, t_s)$ is $c$. Then set
$D_n(c)^-=zD_n(c)^+$. Thus $D_{n}(c)^{\pm}$
are the associated split conjugacy classes of type $c^{(n)}$ (see (\ref{E:speclass})).
Denote by $\sigma_n (c)$ the super class function on $\tG n$ which takes 
value $\pm\frac{n}{\sqrt 2}\zeta_c$ 
on elements in the conjugacy classes $D_{n}(c)^{\pm}$, 
and $0$ elsewhere. For
$\rho = \{ i^{m_i (c)} \} 
\in \mathcal {OP}_n (\G_*)$, 
$\sigma_{\rho}$ =$ 
\prod_{i\in 2\mathbb Z_++1, c \in \G_*} 
\sigma_i (c)^{m_i (c)}$
is the class function of $\tG n$ which takes value
$\pm 2^{-l(\rho)/2}Z_{\rho}$ on the conjugacy classes $D_{\rho}^{\pm}$ and
$0$ elsewhere. Given $\g \in R(\G)$, we
denote by $\sigma_n (\g )$ the class function on $\tG n $ which takes 
value $\pm\frac{n}{\sqrt 2}\g (c) $ on $D_n(c)^{\pm}, c \in \G_*$,
and $0$ elsewhere. 

The following lemma is not difficult to verify.
\begin{lemma}  \label{lem_isom}
  The map $ch$ sends $\sigma_{\rho}$ to $a_{ - \rho} '$.
 In particular, it sends $\sigma_n (c)$ to $a_{ -n} (c)$ in $\tSG$ 
 and $\sigma_n(\g )$ to $a_{ -n} ( \g )$ for 
$n\in 2\mathbb Z+1$.
\end{lemma}

In Sect. \ref{S:chartab}, we will see
that the space $\tSG$ has another distinguished 
basis consisting of generalized Schur Q-functions, which give rise to 
some integral basis in $R^-(\tG n)$.
\subsection{The image of $\chi_n (\g )$ under $ch$}
      
   Recall that we have defined a map from 
$R(\G)$ to $R^-(\tG n )$ (Subsection \ref{S:spinrep}).  

\begin{proposition}  \label{prop_exp}
   For any $\g \in R(\G)$, we have
\begin{equation}
 \sum\limits_{n \ge 0}  2^{\overline n/2}\ch ( \chi_n( \g ) ) z^n 
  = \exp \Biggl( \sum_{ n \ge 1, odd} 
      \frac 2n \, a_{-n}(\g )z^n \Biggr), \label{eq_exp}
\end{equation} 
where $\overline n$ is $0$ or $1$ according to $n$ is even or odd. 
\end{proposition}
\begin{proof} The character value of $\chi_n(\g)$ is given in
Corollary \ref{C:char}, and we have
   \begin{eqnarray*}
  \sum\limits_{n \ge 0}  2^{\overline n/2}\ch ( \chi_n( \g ) ) z^n 
  &= & \sum_{\rho} 2^{l(\rho)}Z_{\rho}^{ -1} 
         \prod_{c\in \G_*} \g (c)^{l (\rho(c))}
          a_{ -\rho (c) }' z^{|| \rho||}                  \\
  &= & \prod_{c\in \G_*} \Bigl ( \sum_{\lambda }
         (2\zeta_c^{ -1}\g (c) )^{l (\lambda)}
         z_{\lambda}^{-1} a_{- \lambda} (c) z^{|\lambda|} \Bigr )    \\
  &= & \exp \Biggl  ( \sum\limits_{ n \geq 1}
         \frac2n \sum\limits_{c \in \G_*}
           \zeta_c^{ -1} \g(c) a_{-n} (c) z^n \Biggl )      \\
  &= & \exp \Biggl( \sum_{ n \ge 1} 
         \frac 2n \, a_{-n}(\g )z^n \Biggr).
 \end{eqnarray*}
 
 Let  $\beta, \g $ be the characters of two representations
 of $\G$. 
 It follows from (\ref{eq_virt}) that 

 \begin{eqnarray*}
  & & \sum\limits_{n \ge 0}  2^{\overline n/2}
\ch ( \chi_n( \beta -\g ) ) z^n  \\
  &=& \Biggl (
       \sum\limits_{n \ge 0} 2^{\overline n/2} \ch ( \chi_n( \beta ) ) z^n
      \Biggl )  \cdot 
      \Biggl (
       \sum\limits_{n \ge 0} 2^{\overline n/2}
 \ch ( \chi_n( \g ) ) ( -z)^n 
      \Biggl )     \\
  &=& \exp \Biggl( \sum_{ n \ge 1, odd} 
         \frac 2n \, a_{-n}(\beta - \g )z^n \Biggr).
 \end{eqnarray*}
Therefore the proposition holds for $\beta - \g$, and so 
for any element $\g \in R_{\mathbb Z}(\G)$.
\end{proof}

\begin{corollary}  \label{cor_char}
   The formula (\ref{E:charvalue}) holds for any $\g \in R(\G)$.
 In particular $\chi_n (\wt)$ is self-dual if $\xi$ is self-dual.
\end{corollary}
 
Component-wise, we obtain 
\begin{equation*}
  \ch (\chi_n (\g ) )
  = 2^{-\overline n/2}\sum\limits_{\rho } \frac {2^{l(\rho)}}{z_\rho }\,
             a_{-\rho}(\g ), 
\end{equation*}
where the sum runs through all the partitions 
$\rho$ of $n$ into odd integers. 
\subsection{Isometry between $\tRG$ and $\tSG$}
    It is well known that there exists a natural Hopf
algebra structure on the symmetric algebra
$\tSG$ with the usual multiplication
and the comultiplication $\Delta$ characterized by
\begin{equation}\label{E:hopf3}
  \Delta ( a_{-n} (\g ))
   = a_{-n} (\g ) \otimes 1 + 1 \otimes a_{-n} (\g ), \qquad n\in 2\mathbb Z_++1.
\end{equation}
Recalling the  Hopf algebra structure
on $\tRG$ defined in Sect.~\ref{S:Hopf}, we can easily
verify the following 
proposition as in the untwisted case.

\begin{proposition}
  The characteristic map $ \ch: \tRG \longrightarrow \tSG$
 is an isomorphism of Hopf algebras.
\end{proposition}
\begin{proof} By counting dimensions of homogeneous degree
subspaces it is easy to see that
$ch$ is an isomorphism of vector spaces. 
The algebra isomorphism
follows simply from the Frobenius reciprocity.
To check the coalgebra isomorphism we use Proposition \ref{prop_exp} to pass from the generators $a_n(\g)$
to the character $\chi_n(\g)$. It is then
a simple calculation to verify that $\chi_n(\g)$
is group-like under the comultiplication (\ref{E:hopf2}), and this shows that $ch$ is a Hopf algebra isomorphism by using
(\ref{E:hopf3}).
\end{proof}

Recall that we have defined a bilinear form 
$\langle \  ,  \, \rangle_{\wt }$ on $\tRG$ and
a bilinear form on $\tSG$ denoted by
$\langle \ , \, \rangle_{\wt }'$. The 
lemma below follows from our definition of
$\langle \ , \, \rangle_{\wt }'$ and the
comultiplication $\Delta$.

\begin{lemma}
  The bilinear form $\langle \  ,  \, \rangle_{\wt } '$ on $\tSG$
 can be characterized by the following two properties:
 
 1). $\langle a_{ -n} (\beta ), a_{ -m} (\g ) \rangle_{\wt}'
  = \frac{n}{2}\delta_{n, m} \langle \beta , \g  \rangle_{\wt}' ,
  \quad \beta, \g \in \G^*, m, n\in 2\mathbb Z_++1.$

 2). $ \langle f g , h \rangle_{\wt}'
       = \langle f \otimes g, \Delta h \rangle_{\wt}' ,$
 where $f, g, h \in \tSG $, and the bilinear form on the r.h.s
 of 2), which is defined on $\tSG \otimes \tSG$, is induced
 from $\langle \ , \ \rangle_{\wt}'$ on $\tSG$. 
\end{lemma}

\begin{theorem}  \label{th_isometry}
  The characteristic map $\ch$ is an isometry from the space
 $ (\tRG, \langle \ \ , \ \  \rangle_{\wt } )$ to
 $ (\tSG, \langle \ \ , \ \  \rangle_{\wt }' )$.
\end{theorem}

\begin{proof} Let $f$ and $g$ be any two super
class functions in $R^-(\tG n)$.
By definition of ch (\ref{E:ch}) it follows that
\begin{align*}
&\langle ch(f), ch(g)\rangle_{\xi}'\\
&=\sum_{\rho, \rho'\in \mathcal{OP}_n (\G_*)}
\frac{2^{(l(\rho)+l(\rho'))/2}}{Z_{\rho}Z_{\rho'}}
f(\rho)g(\rho')\langle a_{-\overline{\rho}}', a_{-\overline{\rho'}}'\rangle_{\xi}'\\
&=\sum_{\rho, \rho' \in \mathcal{OP}_n (\G_*)}
\frac{2^{(l(\rho)+l(\rho'))/2}}{Z_{\rho}Z_{\rho'}}
f(\rho)g({\rho}')\frac{Z_{\rho'}}{2^{l(\rho')}}
\prod_{c\in\G_*}\xi(c)^{l(\rho(c))}
\delta_{\rho, \overline{\rho'}}\\
&=\sum_{\rho\in \mathcal{OP}_n (\G_*)}
\frac1{Z_{\rho}}f(\rho)g(\overline{\rho})
\prod_{c\in\G_*}\xi(c)^{l(\rho(c))}\\
&=\langle f, g\rangle_{\xi},
\end{align*}
where we have used the inner product identity
(\ref{eq_inner}).
\end{proof}
\begin{remark} We can also prove it by showing that
the characteristic map preserves the inner product
of basis elements $\sigma_{\rho}\in \tRG$ and that of
$a_{-\rho}=ch(\sigma_{\rho})\in \tSG$ as in \cite{FJW1}. 
\end{remark}

From now on we will identify the inner product
$\langle \ , \ \rangle_{\wt}$ on $\tRG$
with the inner product $\langle \ , \ \rangle_{\wt} '$ on $\tSG$.
As a special case,  the standard Hermitian form on $R^-(\G_n)$ and therefore
on $\tRG $ is compatible via the characteristic map $\ch$
with the Hermitian form characterized by (\ref{eq_bili})
on $\tSG$.

\section{Vertex operators and $\tRG$}
\label{sect_vertex}
\subsection{A central extension of $\Rz/2\Rz$}
    From now on we assume that $\xi$
is a self-adjoint virtual character
of $\G$, and thus $\Rz$ is an integral lattice
under the symmetric bilinear form $\langle \ , \ \rangle_{\xi}$.

    Let $2\Rz$ be the sublattice of $\Rz$ consisting of elements $2\alpha, \alpha\in\Rz$. The quotient $\Rtz=\Rz/2\Rz$ has an induced abelian group structure
and it can also be viewed
as an $(r+1)$-dimensional vector space over 
$\mathbb F_2=
\mathbb Z/2\mathbb Z$. We will denote by $\overline{\alpha}$
the natural image of $\alpha$ in $\Rtz$.
Define $c_1$ to be the alternating form: 
$\Rtz\times \Rtz\rightarrow \mathbb F_2$ given by
$c_1(\oa, \ob)=\langle \alpha ,\be \rangle_{\xi}+\langle\alpha ,\alpha\rangle_{\xi}
\langle\be ,\be\rangle_{\xi}\, (mod \, 2)$, and let $r_0$ be its rank over $\mathbb F_2$.

The alternating form $c_1$ gives rise to a central 
extension 
$\Rtzh$ of the abelian group $\Rtz$ by the two-element
group $\langle
\pm 1\rangle$ (see \cite{FLM1}):
\begin{equation}
1\rightarrow <\pm 1>\hookrightarrow \Rtzh
\stackrel{\Breve{}}{\rightarrow}\Rtz
\rightarrow 1,
\end{equation}
such that $aba^{-1}b^{-1}=(-1)^{c_1(\breve{a}, \breve{b})}$, $a, b\in \Rtzh$. 

The elements of $\Rtzh$ can be presented as
$\pm e_{\overline{\alpha}}$, where
$\alpha\in \Rz$, which implies that $dim(\Rtzh)=2^{r+2}$.
We note that $e_{\oa}\in \Rtzh$ satisfies
$(e_{\oa})^2=1$. 

Let $\Phi$ be a subgroup of $\Rz$ which is maximal
such that the alternating form $c_1$ vanishes on $\Phi/2\Rz$. A
variant of the following lemma was given in \cite{FLM2}.

\begin{lemma}\label{L:centralext}
There are $2^{(r+1-r_0)}$ irreducible
$\Rtzh$-module structures on the
space $\mathbb C[\Rz/\Phi]$ such that $-1\in\Rtzh$ acts faithfully and 
\begin{equation} \label{E:centralext}
e_{\oa}e_{\ob}=e_{\ob}e_{\oa}(-1)^{c_1(\oa, \ob)}
\end{equation}
as operators on $\mathbb C[\Rz/\Phi]$.
The dimension of $\mathbb C[\Rz/\Phi]$ is equal to $2^{\frac 12r_0}$.
\end{lemma}
We will denote the elements of $\mathbb C[\Rz/\Phi]$
by $e^{[\alpha]}$,
where $[\alpha]=\alpha+\Phi \in \Rz/\Phi$. Clearly
\begin{equation*}
e^{2[\alpha]}=1, e^{[\alpha+\beta]}=e^{[\alpha]}e^{[\beta]}.
\end{equation*}

For $\alpha, \beta\in \Rz$ we write the action
of $\Rtzh$ on $\mathbb C[\Rz/\Phi]$ as
\begin{equation}\label{E:cocycle2}
e_{\oa}.e^{[\beta]}=\ep(\alpha, \beta)e^{[\alpha+\beta]}.
\end{equation}
Then one can check that $\ep$ is a well-defined
cocycle map 
 from $\Rz\times \Rz\rightarrow \langle\pm 1\rangle$. 
One also has $\ep(\alpha, \beta)=\ep(\alpha, -\beta)$. 
\subsection{Twisted Vertex Operators $X ( \g, z)$}

Fix an irreducible $\Rtzh$-module structure on
$\mathbb C[\Rz/\Phi]$
described in Eqn. (\ref{E:cocycle2}). 

We extend the actions of $e_{\oa}$ to the space of tensor product
$$\tFG = \tRG \bigotimes \mathbb C[\Rz/\Phi], 
$$ 
by letting them act on the $\tRG$ part trivially.

Introduce the operators $ H_{ \pm n}( \g ), 
\g \in R(\G), n > 0 $ 
as the following compositions of maps:
\begin{eqnarray*}
  H_{ -n} ( \g ) &:&
    R^- ( \tG m ) 
  \stackrel{ 2^{\overline n/2}\chi_n (\g) \otimes}{\longrightarrow}
    R^- ( \tG n\htimes\tG m )
  \stackrel{ {Ind} }{\longrightarrow}
    R^- ( \tG {n +m} )   \\
  H_n ( \g ) &:&
    R^- ( \tG m ) 
   \stackrel{ {Res} }{\longrightarrow}
    R^-( \tG n\htimes\tG {m -n})
   \stackrel{ \langle 2^{\overline n/2}\chi_n (\g), 
\cdot \rangle_{\wt} }{\longrightarrow}
    R^-( \tG{m -n}) .
\end{eqnarray*}

Define 
\begin{equation*}
  H_+ (\g, z) = \sum_{ n > 0} H_{ -n} ( \g ) z^n, \quad
  H_- (\g, z) = \sum_{ n > 0} H_{ n} (\g )z^{ -n} .
\end{equation*}
  We now define the twisted vertex operators 
$X_n (\g ), n \in {\mathbb Z}$, $\g\in\RG$
by the following generating functions:
\begin{eqnarray}  \label{eq_vo}
 X^+ ( \g, z)
 & \equiv& X ( \g, z) \\
  & =& \sum\limits_{n \in 
  {\mathbb Z}}
 X_n( \gamma) z^{ -n}  \nonumber     \\
  & =&  H_+ (\g , z) H_- (\g , -z) e_{\overline{\g}}.
\nonumber       
\end{eqnarray}

We also denote
\begin{eqnarray*}
 X^- ( \g, z) 
  & \equiv & X ( -\g, z)=X(\ga, -z)  \\
  & =& \sum\limits_{n \in \mathbb Z}
  X^-_n( \g )z^{-n}.  
\end{eqnarray*}
The operators $X_n (\g )$
are well-defined operators acting on the space $ \tFG.$
We extend the bilinear form 
$\langle \ , \ \rangle_{\wt}$ on $\tRG$ to $\tFG $ by letting 
\[
  \langle f e^{[\alpha]}, g e^{[\beta]}\rangle_{ \wt} =
   \langle f, g \rangle_{\wt} \delta_{[\alpha],[\beta]}, \quad
    f, g\in \tRG, \alpha, \beta\in\Rz.
\]
We extend the $\mathbb Z_+$-gradation 
from $\tRG$ to 
$\tFG$ by letting
\begin{eqnarray*}
  \deg  a_{ -n} (\g ) = n , \quad
  \deg e_{\overline{\g} } = 0.
\end{eqnarray*}

Similarly we extend the bilinear form 
$\langle \  , \ \rangle_{\wt }$  to the space
$$
  \tVG = \tSG \bigotimes \mathbb C[\Rz/\Phi]
$$
and extend the $\mathbb Z_+$-gradation on $\SG$
to a $\mathbb Z_+$-gradation on $\tVG$.

The characteristic map $\ch$ will be extended
to an isometry from $\tFG$ to $\tVG$ by 
fixing the subspace $\mathbb C[\Rz/\Phi]$.
We will denote this map again by $\ch$.
\subsection{Twisted Heisenberg algebra and $\tRG$}
    We define $ \widetilde{a}_{ -n} (\gamma), n\in 2\mathbb Z_++1
$ to be a map 
from $\tRG$ to itself by the following composition
\[  R^- (\tG m) \stackrel{ \sigma_n ( \g ) \otimes }{\longrightarrow}
  R^-(\tG n) \bigotimes R^- (\tG m)  \stackrel{{Ind} }{\longrightarrow}
  R^- ( \tG{n +m}).
\]
We also define $ \widetilde{a}_{ n} (\gamma), n\in 2\mathbb Z_++1$ to be a map from $\tRG$ 
to itself
as the composition
\[
  R^-(\tG m)  \stackrel{ Res }{\longrightarrow}
   R^-(\tG n)\bigotimes R^- ( \tG{m -n})
 \stackrel{ \langle \sigma_n ( \g), \cdot \rangle_{\wt}}{\longrightarrow}
 R^- ( \tG {m -n}).
\]
We denote by $\RGO $ the radical of the bilinear form
$\langle \ , \ \rangle_{\wt}$ in $\tRG$ and denote by
$\RGG$ the quotient $\tRG / \RGO$, which 
inherits the
bilinear form $\langle \ , \ \rangle_{\wt}$ from $\tRG$.
          
\begin{theorem}  \label{th_heis}
  $\tRG$ is a representation of the twisted Heisenberg algebra 
 $\thg $  by letting $ a_n (\g )$ 
 $( n \in 2\mathbb Z+1)$ act as $ \widetilde{a}_{ n} (\g )$
and $C$ as $1$. $\RGO$ is a subrepresentation
 of $\tRG$ over $\thg$ and the quotient $\RGG $
 is irreducible. The characteristic map 
 $\ch$ is an isomorphism of $\tRG$ (resp. $\RGO$, $\RGG$)
 and $\tSG$ (resp. $\SGO$, $\tSGG$) as representations over $\thg$. 
\end{theorem}

\subsection{The characteristic map of twisted vertex operators}
    We extend the characteristic map $ch$ to a linear map
$ch$: $End(\tRG)\rightarrow End(\tSG)$ by
\begin{equation}
ch(f).ch(v)=ch(f.v), f\in End(\tRG), v\in\tRG.
\end{equation}
The relation between the vertex operators defined 
in (\ref{eq_vo}) and the Heisenberg algebra $\hg $ is revealed
in the following theorem.

\begin{theorem} For any $\g \in R(\G)$, we have
  \begin{eqnarray*}
   \ch \bigl ( H_+ (\g, z) \bigl )
   &=& \exp \biggl ( \sum\limits_{ n \ge 1, \, odd} \frac 2n \,
    a_{-n} ( \g ) z^n \biggr ), \\
   \ch \bigl ( H_- (\g , z) \bigl )
   &=& \exp \biggl ( \sum\limits_{n \ge 1,\, odd}\frac 2n \,
    a_n (\g ) z^{-n}\biggr ).
  \end{eqnarray*}
\end{theorem}

\begin{proof} Observe that the operator  
$H_+ (\g , z)$
 is the adjoint operator of $ H_- (\g , z^{ -1})$ 
with respect to the bilinear form
 $\langle \ , \ \rangle_{\wt}$.
 Then the theorem follows from 
Lemma~\ref{lem_isom} and  Proposition~\ref{prop_exp} 
by invoking the characteristic map.
\end{proof}

As a consequence we have
\begin{eqnarray*}     
  && \ch \bigl ( X( \g , z)\bigl )  \\
  &= & \exp \biggl ( \sum\limits_{ n \ge 1, \, odd} 
  \frac 2n \, a_{-n} ( \g ) z^n \biggr ) \,
  \exp \biggl ( -\sum\limits_{ n \ge 1, \, odd}
  \frac 2n \,{ a_n( \g)} z^{ -n} \biggr )e_{\overline{\g}}.     
\end{eqnarray*} 
Thus the characteristic map identifies the twisted vertex operators
$X(\g, z)$ defined via finite groups $\tG n$ with
the usual twisted vertex operators of \cite{FLM1, FLM2}.   
\section{Vertex representations and the McKay correspondence}
\label{sect_ade}
\subsection{Product of two vertex operators}
The normal ordered product $: X(\alpha, z) X(\beta, w) :$,
$\alpha, \beta \in R(\G)$ of two vertex operators is defined
as follows:
\begin{equation*}
: X(\alpha, z) X(\beta, w) : =
H_+(\alpha, z)H_+(\beta, w) H_-(\alpha, -z)H_-(\beta, -w)
e_{\oa+\ob}.
\end{equation*}

 In the following theorem and later 
the expression $\big(\frac{z -w}{z+w}
\big)^{ \langle \alpha, \beta \rangle_{\wt}}$
represents the power series expansion in the variable
$\frac wz$. 
\begin{theorem}  \label{th_ope}
 For $\alpha, \beta \in R(\G)$ one has the following
operator product expansion identity for
twisted vertex operators.
 \begin{eqnarray*}
  X(\alpha, z) X(\beta, w) & =& \ep(\alpha, \beta)
  :X(\alpha, z) X(\beta, w):
  \biggr(\frac{z -w}{z+w}\biggr)^{ \langle \alpha, \beta \rangle_{\wt}}. 
 \end{eqnarray*}
\end{theorem}
\begin{proof}It suffices to compute that
 \begin{eqnarray*}
   && ch(H_- (\alpha, -z) H_+ (\beta, w))    \nonumber \\
   & =& \exp\,\biggl ( -\sum\limits_{ n \ge 1, \, odd}
    \frac 2n \,{ a_n(\alpha)} z^{ -n} \biggr )
   \exp \biggl ( \sum\limits_{ n \ge 1, \, odd} \frac 2n \,
    a_{-n} (\beta) w^n \biggr )    \nonumber  \\ 
  & =& ch(H_+ (\beta, w)H_- (\alpha, -z) )  
 \exp \, \biggl  ( - \langle \alpha, \beta \rangle_{\wt}
       \sum\limits_{ n \ge 1, \, odd} \frac 2n \,
       z^{ -n} w^n \biggr )  \nonumber  \\ 
  & =&  ch(H_+ (\beta, w)H_- (\alpha, -z) ) 
        \biggr(\frac{z-w}{z+w}\biggr)^{\langle \alpha, \beta \rangle_{\wt} } .
      \nonumber 
  \end{eqnarray*}
\end{proof}
 
The following proposition is easy to check.

\begin{proposition}  \label{prop_prim}
  Given $\alpha \in R(\G), \beta \in \Rz$ and $n\in 2\mathbb Z+1$, we have
  $$ 
   [ a_n (\alpha ), X(\beta, z)]
   = \langle \alpha, \beta \rangle_{\wt} X(\beta, z) z^n.
  $$
\end{proposition}   
\subsection{Twisted affine Lie algebra $\tloopg$ 
and twisted toroidal Lie algebra $\thhg$}
  Let $\mathfrak g$ be a rank $r$ complex simple Lie algebra of ADE type, and
let $\overline{\Delta}$ be the root system generated
by a set of simple roots 
$\alpha_1. \ldots, \alpha_r$. Let $\alpha_{max}$ be the highest root.
The Lie algebra is generated by
the Chevalley generators $e_{\alpha_i}, e_{-\alpha_i}, 
h_i=h_{\alpha_i} $.
We normalize the invariant bilinear form on  $\mathfrak g$ 
by $(\alpha_{max},
\alpha_{max}) = 2$. 

Let $\theta$ be an automorphism of $\mathfrak g$ of order $k$
and let $\omega=exp(2\pi i/k)$. The automorphism 
$\theta$ induces a 
$\mathbb Z/k\mathbb Z$-gradation for $\mathfrak g$:
\begin{equation*}
\gs=\bigoplus_{i\in\mathbb Z/k\mathbb Z}\gs_i, \qquad
\gs_i=\{g\in \gs\vert \theta(g)=\omega^i g\},
\end{equation*} 
and
\begin{equation*}
[\gs_i, \gs_j]\subset \gs_{i+j}.
\end{equation*}

The twisted affine Lie algebra $\widehat{\gs}[\theta]$
is the graded vector space
\begin{equation}
\widehat{\gs}[\theta]=\bigoplus_{i=0}^{k-1}\gs_i \otimes 
t^i\mathbb C[t^k, t^{-k}]\bigoplus \mathbb CC
\end{equation}
with the commutating relations
\begin{eqnarray} \label{E:twrel}
  [ a(n), b(m)] & =& [a, b] (n+m) + \frac nk \delta_{n, -m} (a| b) C , \\
  {[C, a(n)]}   & =& 0, \quad a, b \in \mathfrak g, n, m \in \mathbb Z,
\end{eqnarray}
where we used the notation 
$$
 a (n) = a \otimes t^n, \quad a \in \mathfrak g, n \in \mathbb Z.
$$
 When $\theta$ is the identity, $\loopg[1]$ becomes the (untwisted)
affine Lie algebra $\loopg$.
Let $A = (a_{ij})_{ 0 \leq i,j \leq r}$ be the
affine Cartan matrix associated to $\loopg$.
The submatrix $(a_{ij})_{ 1 \leq i,j \leq r}$ is the Cartan matrix 
of $\mathfrak g$.

The linear map
\begin{align*}
e_{\alpha_i}&\longrightarrow e_{-\alpha_i}\\
h_{\alpha_i}&\longrightarrow -h_{\alpha_i}
\end{align*}
defines an involution of the Lie algebra $\mathfrak g$.
We will denote the associated twisted affine Lie algebra
by $\tloopg$. Let
$k_{\alpha}=\frac 12(e_{\alpha}+e_{-\alpha})$
and $p_{\alpha}=\frac 12(e_{\alpha}-e_{-\alpha})$, where
$\alpha\in\overline{\Delta}$. It is easily
seen that 
\begin{align*}
\mathfrak g_0&=\bigoplus \mathbb Ck_{\alpha},\\
\mathfrak g_1&=\bigoplus \mathbb Cp_{\alpha}\oplus
\bigoplus \mathbb Ch_{\alpha}
\end{align*}

The {\em basic twisted representation}  $V$ of $\tloopg$ is
the irreducible highest weight representation generated by
a highest weight vector which is annihilated by 
$ a(n), n \in \mathbb Z_+, a \in \mathfrak g$ and
$C$ acts on $V$ as the
identity operator. 

We now introduce the complex twisted toroidal Lie algebra
$\thhg$ 
(associated to $\mathfrak g$) with the
following presentation: the
generators are 
$$
 C, h_i (m), x_n (\pm\alpha_i), m\in 2\mathbb Z+1,
n \in \mathbb Z, i =0, \ldots, r;
$$
and the relations are given by:
$C$ is central, and
\begin{align}
   &x_n(\alpha_i)= (-1)^n x_n(-\alpha_i), \nonumber\\ 
   &{[ h_i (m), h_j (m')]} 
  =\frac m2 a_{ij} \delta_{m, -m'}C,    \nonumber  \\
   &{[h_i (n), x_m (\alpha_j)]}
  =a_{ij} x_{n +m} (\alpha_j),  \nonumber    \\
   &{[x_n (\alpha_i), x_{n'}(- \alpha_i)]}
  = 8\{ h_i (n +n') + n\delta_{n, -n'} C \},
      \label{eq_pres}  \\
&{\sum_{s=0}^{a_{ij}}\binom{a_{ij}}{s}[x_{n+s}(\alpha_i), x_{n'-a_{ij}-s}(\alpha_j)]}=0, 
\quad\mbox{if $a_{ij}\geq 0$} \nonumber \\
&{\sum_{s=0}^{-a_{ij}}(-1)^s
\binom{-a_{ij}}{s}[x_{n+s}(\alpha_i), x_{n'-a_{ij}-s}(\alpha_j)]}= 0, 
\quad\mbox{if $a_{ij}<0$} \nonumber   
\end{align}
where $n, n' \in \mathbb Z$, $m, m'\in 2\mathbb Z+1$, $i, j = 0, 1, \ldots, r$, and 
$h_i(2n)=0$ for $n\in\mathbb Z$. The twisted toroidal algebra is the $q\rightarrow 1$ limit
of the twisted quantum current algebra \cite{J2} (see
a slightly different form in \cite{DI}).

Set $k_i(2n)=\frac 14x_{2n}(\alpha_i)$ and 
$p_i(2n+1)=\frac 14x_{2n+1}(\alpha_i)$. One can check that
the relations given in (\ref{eq_pres}) are consequences of the twisted
algebra 
$\tloopg$ for the case of 
$\theta=-1$ defined in (\ref{E:twrel}) (cf. the proof 
of Theorem \ref{th_mck} later) .

%
%
%
%
%
\subsection{Realization of twisted vertex representations}
Let $\G$ be a finite
subgroup of $SL_2(\mathbb C)$ and the virtual character $\wt$ 
to be twice the trivial character minus the character of the
two-dimensional defining representation of $\G\hookrightarrow SL_2(\mathbb C)$.
It is known \cite{Mc} that the matrix $A = (a_{ij})_{0 \leq i,j \leq r}$
in Sect.~\ref{subsec_mcka}
is the Cartan matrix for the corresponding affine Lie
algebra $\loopg$. 
 
The following theorem provides a finite group
realization of the vertex representation of
the twisted toroidal Lie algebra $\thhg$ on $\tFG$.

\begin{theorem}  \label{th_mck}
 A vertex representation of the twisted toroidal Lie algebra $\thhg$
 is defined on the space $\tFG$ by letting
 \begin{align*}
  x_n (\alpha_i) \mapsto X_n (\g_i), &\qquad 
  x_n (-\alpha_i) \mapsto  \ep(\g_i, \g_i)X_n (-\g_i),\\
  h_i (m) \mapsto a_m (\g_i), &\qquad C\mapsto 1,
 \end{align*}
 where $n \in \mathbb Z, m\in 2\mathbb Z+1, 0 \leq i \leq r$. 
\end{theorem}

\begin{proof} All the commutation relations without binomial
coefficients are easy consequences
of Proposition~\ref{prop_prim} and Theorem~\ref{th_ope} 
by the usual vertex operator calculus in the twisted picture
 (see \cite{FLM2, J}). The corresponding relations with binomial coefficients
in $\tVG$ are equivalent to 
\begin{align*}
(z+w)^{a_{ij}}[X(\g_i, z), X(\g_j, w)]&=0, \qquad \mbox{$a_{ij}\geq 0$},\\
(z-w)^{-a_{ij}}[X(\g_i, z), X(\g_j, w)]&=0, \qquad \mbox{$a_{ij}<0$}.
\end{align*}
This is again proved by using Theorem~\ref{th_ope} with the same method
as in the quantum vertex operators \cite{J2}.
\end{proof}

 Recall that $\delta = \sum_{i =0}^r d_i \g_i$ generates the 
one-dimensional radical $R^0_{\mathbb Z}$ of the bilinear form
$\langle \ , \ \rangle_{\wt}$ in $\Rz$, where
$d_i$ is the degree of the irreducible character $\g_i$ of $\G$. 
The lattice $\Rz$ in this case
can be identified with the root lattice
for the corresponding affine Lie algebra. 
The quotient lattice $\Rz / R^0_{\mathbb Z}$ inherits 
a positive definite integral bilinear form. 
Denote by $\Gbar$ the set of non-trivial irreducible characters of $\G$:
$$ 
  \Gbar = \{ \g_1, \g_2, \ldots, \g_r \}.
$$
Let $R_{\mathbb Z}(\Gbar)$ be the sublattice of $\Rz$ generated by
$\Gbar$. 
Denote by $Sym(\Gbar)$ the symmetric algebra generated by
$a_{-n} (\g_i),$ $ n\in 2\mathbb Z_++1, $ $i =1, \ldots , r$. Equipped with the bilinear
form $\langle \ , \ \rangle_{\wt}$, $Sym (\Gbar)$
is isometric to $\tSGG$ which is in turn isometric
to $\RGG$ as well. The irreducible $\Rtzh$-module $\mathbb C[\Rz/\Phi]$
induces an irreducible $\Rtzhbar$-module structure on
$\mathbb C[R_{\mathbb Z}(\Gbar)/\overline{\Phi}]$ given by
the restriction of the alternating form $c_1$. We let $\overline{r}_0$
denote the rank of the restriction of $c_1$, then the statement of 
Lemma \ref{L:centralext} also holds for the sublattice $\Rtzbar$ and $\Rtzhbar$.
In this case if the determinant of the Cartan matrix is an odd integer, then
$\overline{r}_0=0$ and the space $\mathbb C[R_{\mathbb Z}(\Gbar)/\overline{\Phi}]$
is trivial. 
 
We define
\begin{align*}
 \tVGG  &= \tSGG \bigotimes \mathbb C[R_{\mathbb Z}(\Gbar)/\overline{\Phi}] 
 \cong Sym (\Gbar) \bigotimes \mathbb C[R_{\mathbb Z}(\Gbar)/\overline{\Phi}] ,\\
\tFGG  &= \tRGG \bigotimes \mathbb C[R_{\mathbb Z}(\Gbar)/\overline{\Phi}].
\end{align*}
Obviously $ch$, when restricted to $\tFGG$, is an isometric
isomorphism onto 
$\tVGG$.

 The space $\tFG$ associated to the lattice
$\Rz$ is isomorphic to the tensor
product of the space $\FGG$ 
associated to $\Rzz$ and the space associated to
the rank $1$ lattice $\mathbb Z \delta$
equipped with the zero bilinear form.

The identity for a product of vertex operators $X (\g, z)$ 
associated to $\g \in \overline{\Delta}$ 
(cf.~Theorem~\ref{th_mck}) implies that
$\tVGG$ provides a realization of the vertex representation
of $\tloopg$ on $\tVGG$ (cf. \cite{FLM1}).
The following theorem establishes
a direct link from the finite group $\G \in SL_2(\mathbb C)$
to the affine Lie algebra $\tloopg$. This  gives
a twisted version of the new form of the
McKay correspondence given in \cite{FJW1}.

\begin{theorem}
 The operators $X_n(\g ), \g \in \overline{\Delta}, a_n (\g_i),
 i =1, 2, \ldots, r, n \in \mathbb Z$
 define an irreducible representation of the affine Lie algebra $\tloopg$
 on $\tFGG$ isomorphic to the twisted basic representation. 
\end{theorem}
\section{Vertex operators and irreducible characters of $\tG n $}
\label{sect_char}
     In this section we specialize $\wt$ to be 
the trivial character $\g_0$ of $\G$. We will
describe how to obtain the character table for the spin
supermodules of $\tG n$ from our vertex operator
approach, generalizing \cite{J}.
\subsection{Algebra of vertex operators for $\wt = \g_0$}
      In this case the weighted bilinear form reduces 
to the standard one $\langle \ , \ \rangle$
and $\Rz$ is isomorphic to the lattice $\mathbb Z^{r+1}$
with the standard integral bilinear form. Recall that
$\langle \g_i , \g_j \rangle = \delta_{ij}$. 
For simplicity we will only consider the vertex representations on the space $\RG$. 
In the following result
the bracket
$\{ \ , \ \}$ denotes the anti-commutator.

\begin{theorem}  \label{th_cliff}
  The operators $X^+_n (\g_i), X^-_n (\g_i)$
 $(n \in \mathbb Z, 0 \leq i \leq r)$
 generate a generalized Clifford algebra:
 \begin{eqnarray}
  {[X^+_n( \g_i), X^+_{n'}(\g_j)]} &= & 0, \quad\mbox{$i\neq j$} \nonumber\\
{[X^-_n( \g_i), X^-_{n'}(\g_j)]} &= & 0, \quad\mbox{$i\neq j$} \nonumber\\
\{X^+_n( \g_i), X^+_{n'}(\g_i)\} &=&2(-1)^n\delta_{n, -n'}, 
\nonumber   \\
\{X^-_n( \g_i), X^-_{n'}(\g_i)\} &=&2(-1)^n\delta_{n, -n'}, 
\label{eq_delta}   \\
{[X^+_n(\g_i), X^-_{n'}(\g_j)]} &= & 0, \quad\mbox{$i\neq j$} \nonumber\\
  \{X^+_n(\g_i), X^-_{n'}(\g_i)\} &= & 2\delta_{n, -n'}.
  \nonumber
 \end{eqnarray}
\end{theorem}
\begin{proof}
It follows from the standard vertex operator
calculus (cf. \cite{J}) by using Theorem~\ref{th_ope}.
\end{proof}
Therefore we see that $\tRG$ is isomorphic
to the tensor product of $r+1$ copies of the
space $R^-$, the sum of Grothendieck groups of 
spin characters of $\tS n$-supermodules.

\begin{remark} Let $A=(1-\delta_{ij})_{(r+1)\times(r+1)}$, the 
matrix of the alternating
form $c_1$ over $\mathbb F_2$, then $A^2=\overline rI$.
Here $\overline r=0 $ if $r$ is even and $1$ if $r$ is odd.
Consequently it follows from
Lemma \ref{L:centralext}
there are exactly $2^{\overline{r+1}}$ irreducible $\hat{R}^-_{\mathbb F_2}(\G)$-module
structures on the $2^{\lceil\frac{r+1}2\rceil}$-dimensional space 
$\mathbb C[R_{\mathbb Z}(\G)/\Phi]$. 
One of the (at most) two
irreducible module structures is given by the cocycle
$\ep(\g_i, \g_j)=1$,
for $i\leq j$, and $\ep(\g_i, \g_j)=-1$,
for $i>j$. Then the vertex operators 
$X^{\pm}_n(\g_i)$ generate
the twisted Clifford algebra on the space $\tFG$ defined
by $\{X^{\pm}_n(\g_i), X^{\pm}_{n'}(\g_i)\}=2(-1)^n\delta_{ij}\delta_{n, -n'}$ and
$\{X^{+}_n(\g_i), X^{-}_{n'}(\g_i)\}$ $=2\delta_{ij}\delta_{n, -n'}$.
\end{remark}
\subsection{Super spin character tables of $\tG n$ and vertex operators} 
\label{S:chartab}
We now use the twisted vertex operators to construct 
all irreducible characters of spin supermodules of
$\tG n$ for all $n$. 

Let $R^-_{\mathbb Z}(\tG n)$ be the lattice
generated by the characters of spin irreducible 
$\tG n$-supermodules. Then 
$R^-_{\mathbb Z}(\tG n)\otimes\mathbb C\simeq R^-(\tG n)$.

First we construct  a special orthonormal basis in $\tRG$ and then use them 
to give irreducible characters of spin $\tG n$-supermodules. 
The vertex operator $X_n(\g)$ is defined as in (\ref{eq_vo})
except that we drop $e_{\overline{\g}}$.
The following is easily seen (cf. \cite{J}).
\begin{lemma} 
  \label{L:lattice}
 For $n\in \mathbb Z$, $\alpha \in R( \G )$ and
 $ \g \in \G^*$, we have
 \begin{equation*}
  X_{-n}(\pm\g ).1
   =\delta_{n, 0}, \qquad 
    n\geq 0.            
 \end{equation*}
\end{lemma}

For a $m$-tuple index $\phi=(\phi_1, \cdots, \phi_m)\in
\mathbb Z^m$ we denote
\begin{align*}
X_{\phi}(\g )&=X_{\phi_1}(\g )\cdots X_{\phi_m}(\g ).1,\\
x_{\phi}(\g) &=(X_{\phi_1}(\g ).1)\cdots (X_{\phi_m}(\g ).1).
\end{align*}
We also define the raising operator $R_{ij}$ by
\begin{equation*}
R_{ij}(\phi_1, \cdots, \phi_m)=(\phi_1, \cdots, \phi_i+1,
\cdots, \phi_j-1, \cdots, \phi_m).
\end{equation*}
Then we define the action of the raising operator $R_{ij}$
on $X_{\phi}(\g )$ or $x_{\phi}(\g)$ by
$X_{R_{ij}\phi}(\g )$ or $x_{R_{ij}\phi}(\g)$. 

Given $ \la \in \mathcal {OP}(\G^*)$, we define
$$
 X_{\lambda}
  = \prod_{\g \in \mathcal \G^*} X_{-\la( \g )} (\g).
$$
Similarly we define
$x_{\lambda}=\prod_{\g \in \mathcal \G^*}x_{\la( \g )} (\g).$

\begin{theorem}\label{T:char1}
 The vectors $X_{\lambda}$ for
 $ \lambda=(\la(\g ))_{\g \in \G^*} \in \mathcal {SP} (\G^*)$
 form an orthogonal basis in the vector space $\tRG$ with $\langle X_{\lambda}, X_{\mu}\rangle=
2^{l(\lambda)}\delta_{\lambda, \mu}$.
Moreover, we have that
\begin{align} \label{E:raisingop}
X_{\lambda}&=\prod_{\g\in\G_*}
\prod_{i<j}\frac{1-R_{ij}}{1+R_{ij}}X_{\la}(\g)\\
&=x_{\lambda}+\sum_{\lambda\gg\mu}c_{\la, \mu}x_{\mu},
\end{align}
where $c_{\la, \mu}\in \mathbb Z$, and $R_{ij}$ is the raising operator.
\end{theorem}
\begin{proof} The generalized Clifford algebra structure (\ref{eq_delta})
 implies that the non\-zero elements
 $$
  \prod_{\gamma\in \G^*}X_{-n_1}(\gamma)\cdots 
  X_{-n_l}(\gamma).1
 $$
 of distinct indices generate a spanning set in the space $\tRG $.
 To see that they satisfy the raising operator expansion we compute 
 that
 \begin{align*}
   &X (\g , z_1)\cdots X (\g , z_l).1  \\
   & = :X (\g , z_1) \cdots X (\g , z_l): 
      \prod_{i <j} \frac{z_i-z_j}{z_i+z_j} .
 \end{align*}
 Using the result in \cite{J} it follows that this is 
 exactly the generating function of the raising operator expansion 
 at the case $\G$ is trivial under the isomorphism $\ch$. 
In other words, equation (\ref{E:raisingop}) is true when $\la$ is a characteristic partition-valued
function.

Next the orthogonality follows from the generalized Clifford
 algebra commutation relations in Theorem~\ref{th_cliff}.
The orthogonality relations
show that the raising operator is not affected by
the character $\g$, hence the general case follows
by multiplying the raising operator formula
for each $\g$.
\end{proof}

The corresponding basis in $\tSG$ are the classical symmetric
functions called Schur's Q-function \cite{J}.
For any strict partition $\la=(\la_1, \ldots, \la_l)$,
the Schur's Q-function $Q_{\la}$
is determined by \cite{Sc, M} 
\begin{equation}\label{E:schur}
Q_{\la}(y_1,\ldots, y_n)=2^{l}\sum_{\sigma\in
S_n/S_{n-l}}y_{\sigma(1)}^{\la_1}\cdots y_{\sigma(l)}^{\la_l}\prod_{\la_i>\la_j}
\frac{y_{\sigma(i)}+y_{\sigma(j)}}{y_{\sigma(i)}-y_{\sigma(j)}},
\end{equation}
where $S_{n-l}$ acts on $y_{l+1}, \cdots , y_n$
and we allow $\la_j=0$ for $j=l+1, \cdots, n$. 
It is known
(see for example \cite{M, J})
that $Q_{\la}$, $\la\in \mathcal{SP}_n$ form
a basis in the subring of symmetric functions generated by
the power sums $p_1, p_3, p_5, \ldots$.  

We can think of $a_{ -n } (\g ), n >0 , \g \in \G^*$ as
the $n$-th power sum in a sequence of variables 
$ y_{ \g } = ( y_{i\g } )_{i \geq 1}$. In this way
we identify the space $\tSG$ with a distinguished subspace
of 
symmetric functions generated by odd degree power sums indexed by $ \G^*$. 
In particular given a strict partition $\lambda$
we denote by $Q_{\lambda} (\g)$ the Schur's Q-function associated to $y_{\g}$. 
We also denote by $Q_{\lambda}(\g)$ the corresponding element
in $\tSG$ by the identification of $\tSG$ and $\tRG$.
For $\lambda \in \mathcal P (\G^*)$, we denote
\begin{eqnarray}  \label{eq_schur}
  Q_{\lambda} = \prod_{\g \in \G^*} Q_{\lambda(\g)} (\g) \in \tSG.
\end{eqnarray}

For $\la\in\mathcal {SP} ( \G^*)$ we define 
\begin{equation} \label{E:intchar}
\overline{Q}_{\lambda}=
2^{-(l(\la)-d(\la))/2}Q_{\lambda},
\end{equation}
where $d(\la)$ is the parity of $\la$ (see
(\ref{E:parity1}). Similarly we define $\overline{q}_{\la}
=2^{-(l(\la)-d(\la))/2}\prod_{\g \in \mathcal \G^*, i}q_{\lambda_i(\g)}(\g)$, 
where
$q_m(\g)$ = $Q_{(m)}(\g)$.  
In particular we have 
\[
ch(\chi_n(\g))=\overline q_n=2^{-\overline n/2}q_n(\g)
=2^{-\overline n/2}
\overline Q_{(n)}(\g).
\]

The following result immediately follows from Theorem
\ref{T:char1} combined with the characteristic map $ch$ and (\ref{E:intchar}).
\begin{proposition}
For strict partition-valued function $\lambda\in
\mathcal {SP} (\G^*)$ we have
\begin{equation}
\begin{aligned}
\langle \overline Q_{\lambda}, \overline 
Q_{\mu}\rangle &=
0  \qquad\mbox{if $\lambda\neq \mu$},\\
\langle \overline Q_{\lambda}, \overline 
Q_{\lambda}\rangle &=\begin{cases} 1 & 
\mbox{if $\lambda$ is
even,}\\
2 & \mbox{if $\lambda$ is odd.}
\end{cases}
\end{aligned}
\end{equation}
\end{proposition}

The following result generalizes a similar result of
\cite{Jo} for trivial $\G$.
\begin{lemma} Under the characteristic map $ch$
the symmetric functions $\overline q_{\la}$
corresponds to a character in $R^-_{\mathbb Z}(\tG n)$
and $\overline Q_{\la}$ correspond to a virtual character
in $R^-_{\mathbb Z}(\tG n)$.
\end{lemma}
\begin{proof}  Observe that
the tensor product of two irreducible supermodules of type $Q$
is a sum of two irreducible supermodule of type $M$ (see
Proposition \ref{P:supermult}). Computing its 
inner product we know that for positive odd
integers $m$ and $n$ the character 
$\chi_n(\g)\chi_m(\g)$ is twice of some irreducible character.
Then by induction we see that $ch^{-1}(\overline q_{\la})$
is a character in $R^-_{\mathbb Z}(\tG n)$.
From Theorem \ref{T:char1} we see that $\overline{Q}_{\la}$
is a $\mathbb Z$-linear combination of $\overline q_{\mu}$ with
$\lambda\gg \mu$, hence  $\overline{Q}_{\la}$
is a virtual character of $\tG n$.
\end{proof}

For $\la\in\mathcal {SP}_n( \G^*)$, we define
$\tG {\la}=\tG{\la(\ga_0)}\htimes \cdots\htimes \tG{\la(\ga_r)}$.
For a partition $\mu$ and an irreducible character $\ga$ of
$\G$, we define the
spin character $\chi_{\mu}(\g)$ of $\tG{\mu}$ to be 
$\chi_{\mu_1}(\ga)\otimes \cdots\otimes
\chi_{\mu_l}(\ga)$ (see Corollary \ref{C:char}).

\begin{theorem}\label{T:main} 
For each strict partition-valued function $\lambda \in \mathcal {SP}_n( \G^*)$,
the
vector $\overline Q_{\lambda}$ corresponds,
under the characteristic map $ch$,  to the 
irreducible 
character $\chi_{\la}$
of the spin $\tG n$-supermodule given by
a $\mathbb Z$-linear combination of
\begin{equation}\label{E:main3}
Ind_{\tG{\rho}}^{\tG n}
\chi_{\rho(\g_0)}(\g_0)\otimes \cdots\otimes\chi_{\rho(\g_{r})}(\g_r),
\end{equation}
where 
$\rho\leqslant\la$ and the first summand is $\rho=\la$ with multiplicity one.
The parity of $\chi_{\lambda}$ is equal to $d(\la)=
n-l(\la)\, (mod\, 2)$. 
Its character at the conjugacy class of type
$\mu \in \mathcal {OP}_n(\G_*)$ is equal to
 the matrix coefficient
 \begin{eqnarray}\label{E:table}
  2^{(l(\mu)-l(\la)+d(\la))/2}\langle
    X_{\lambda}, a_{-\mu} 
  \rangle. 
 \end{eqnarray}
Moreover the degree of the character 
is equal to
\begin{equation}\label{E:degree}
2^{\lfloor (n-l(\la))/2\rfloor}n!\prod_{\g\in \G^*}
\big(\frac{deg(\g)^{|\la(\g)|}}
{\prod_{1\leq i\leq l(\la(\g))} \la_i(\g)!}
\prod_{i<j}\frac{\la_i(\g)-\la_j(\g)}{\la_i(\g)+\la_j(\g)}
\big),
\end{equation}
where $\lfloor a\rfloor$ denotes the smallest integer $\geq a$.
\end{theorem}
\begin{proof} Suppose we know that
$\overline Q_{\lambda}$ corresponds to the character
of an irreducible $\tG n$-supermodule. By
 (\ref{eq_inner}) and definition of the characteristic
map we see immediately that the linear combination
(\ref{E:main3}) is given by the vertex operator
structure and the matrix coefficient
(\ref{E:table}) give the character table of all
irreducible supermodules.

First we observe that the number of irreducible spin supermodules of
type $M$ is equal to the number of
even strict partition-valued functions, which are realized
by the vectors $2^{-(l(\la)-d(\la))/2}X_{\la}(\g)$ 
($\la\in\mathcal{SP}^0(\G^*))$ up to signs.
As for the vectors $2^{-(l(\la)-d(\la))/2}X_{\la}(\g)$ with
($\la\in\mathcal{SP}^1(\G^*))$, 
it follows from Theorem \ref{T:char1}
and Proposition \ref{P:criterion}
that each of such vectors
corresponds to a virtual irreducible character in 
$R^-_{\mathbb Z}(\tG n)$ of type $Q$, since the case of
sum or difference of two irreducible characters of type M is ruled
out by the orthogonality. To show that they correspond to actually
irreducible characters it is sufficient to
show
that the value of 
$$ch^{-1}(2^{-(l(\la)-d(\la))/2}X_{\la}(\g))$$
at the conjugacy class of the identity element of
$\tG n$ is positive.

Let $c^0\in\G_*$ be the class consisting of the identity
in $\G$. The type of the
identity element in $\tG n$ is the partition-valued function
$\rho$ such that $\rho(c^0)=(1^n)$
and $\rho(c)=0$ for $c\neq c^0$.  
Recall from (\ref{E:classba}) that
$a_m(c^0)=\sum_{\g\in\G^*}deg(\g)a_m(\g)$.
 By comparing weights and using 
orthogonality (Theorem \ref{T:char1})
we have
\begin{align*}
\langle
    X_{\lambda}, a_{-\rho}'
  \rangle &=\langle
    X_{\lambda(\g)}, a_{-1}^{n}(c^0) 
  \rangle\\
&=\langle
    \prod_{\g\in\G^*}X_{\lambda(\g)}(\g), 
\big(\sum_{\g\in\G^*}(deg\g)a_{-1}(\g)\big)^n 
  \rangle\\
&=n!\prod_{\g\in\G^*}\frac{(deg\g)^{|\la(\g)|}}{|\la(\g)|!}
\prod_{\g\in\G^*}\langle
    X_{\lambda(\g)}(\g), a_{-1}^{|\la(\g)|}(\g) 
  \rangle.
\end{align*}

By the result of (6.51) in \cite{J} 
we have
\begin{equation*}
\langle X_{\lambda(\g)}(\g), a_{-1}^{|\la(\g)|}(\g)
\rangle
=\frac{|\la(\g)|!}{\la_1(\g)!\cdots \la_{l(\la(\g))}(\g)!}
\prod_{i<j}\frac{\la_i(\g)-\la_j(\g)}{\la_i(\g)+\la_j(\g)},
\end{equation*}
which implies
 the formula (\ref{E:degree}), thus the theorem is proved. 
\end{proof}

The irreducible spin $\tG n$-supermodules can be described easily
as follows.  For each irreducible character $\g\in\G^*$
let $U_{\g}$ be the irreducible
$\G$-module affording $\g$. For each strict partition $\nu$ let $V_{\nu}$
be the corresponding irreducible spin supermodule of $\tS n$.
Using the construction of Sect. \ref{S:spinrep}, we see that
$U^{\otimes n}_{\g}\otimes V_{\nu}$ is a spin $\tG n$-supermodule.

\begin{proposition} For each strict partition-valued function
$\la=(\la(\g))$ $\in$ 
$\mathcal{SP}_n(\G^*)$, with $m$ of the partitions $\la(\g)$ being odd, 
the 
super tensor product
$$
\prod_{\g\in\G^*} \big(U^{\otimes l(\la(\g))}_{\g}\otimes V_{\la(\g)}\big)
$$ 
decomposes completely
into $2^{\lceil m/2\rceil}$ copies of an 
irreducible spin $\tG{\la}$-super\-module. Denote this
irreducible module by $W_{\la}$. 
Then the induced
supermodule $Ind_{\tG{\la}}^{\tG n}W_{\la}$ 
is the
irreducible spin $\tG n$-supermodule corresponding to $\la$, and it
is of type $M$ or $Q$ according to $d(\la)=n-l(\la)$ is even or odd.
\end{proposition}
\begin{proof} Let $V_{\la}$ be the irreducible spin $\tG n$-supermodule
corresponding to $\la$. It follows from Theorem \ref{T:main} that
$V_{\la}$ is an irreducible component of $Ind_{\tG{\la}}^{\tG n}W_{\la}$.

Note that the supermodule
$U^{\otimes l(\la(\g))}_{\g}\otimes V_{\la(\g)}$
is irreducible and is of type $M$ (or type $Q$) according
to $d(\la(\g))=|\la(\g)|-l(\la(\g))$ even (or odd). 
Let $\la(\g_{i_0})$, $\cdots$, $\la(\g_{i_{m-1}})$
be odd partitions, and let $\la(\g_{i_m})$, $\cdots$, $\la(\g_{i_{r}})$
be even partitions. Then the parity of $\la$ equals the parity
of $m$, i.e., $d(\la)=n-l(\la)\equiv m \, (mod\, 2)$.

It follows from Proposition \ref{P:supermult} 
that $\prod_{\g\in\G^*} U^{\otimes l(\la(\g))}_{\g}\otimes V_{\la(\g)}$ 
decomposes completely into
$2^{\lceil m/2\rceil}$ copies of 
the irreducible supermodule $W_{\la}$ of $\tG{\la}$, and $W_{\la}$ is of type
$M$ if $m$ is even and of type $Q$ otherwise. The degree of
$Ind_{\tG{\la}}^{\tG n}W_{\la}$ equals to $\frac{|\tG n|}{|\tG{\la}|}
deg(W_{\la})$, and we have
\begin{align*}
&deg(W_{\la})=2^{-\lceil m/2\rceil}
\prod_{\g\in\G^*}deg(\g)^{l(\la(\g))}deg(V_{\la(\g)})\\
&=2^{-\lceil m/2\rceil}\prod_{\g\in\G^*}
\big(2^{\lfloor \frac{d(\la(\g))}2\rfloor}
\frac{deg(\g)^{l(\la(\g))}|\la(\g))|!}
{\prod_{1\leq i\leq l(\la(\g))} \la_i(\g)!}
\prod_{i<j}\frac{\la_i(\g)-\la_j(\g)}{\la_i(\g)+\la_j(\g)}\big),
\end{align*}
where we have used the degree formula (\ref{E:degree}) for the
special case of $\tS n$.
The exponents of $2$ in the first factor and the product sum up to 
\begin{align*}
&\sum_{\g\in\G^*}\lfloor d(\la(\g))\rfloor-\lceil m/2\rceil\\
&=\sum_{s=0}^{m-1}\frac{d(\la(\g_{i_s}))+1}2+
\sum_{s=m}^{r}\frac{d(\la(\g_{i_s}))}2-\frac{m-\overline{m}}2\\
&=\frac{n-l(\la)}2+\frac{\overline{m}}2=\lfloor\frac{n-l(\la)}2\rfloor,
\end{align*}
where $\overline m=0$ or $1$ according to $m$ is even or odd.
Thus the degree of $Ind_{\tG{\la}}^{\tG n}W_{\la}$ is exactly
the one given by Eqn. (\ref{E:degree}). Hence
$Ind_{\tG{\la}}^{\tG n}W_{\la}$ is the irreducible spin 
$\tG n$-supermodule $V_{\la}$
corresponding to $\la$.
\end{proof}

\bibliographystyle{amsalpha}

\end{document}